\documentclass{article}
\usepackage{amsmath,amssymb}
\usepackage{cite}
\usepackage[bookmarks=false]{hyperref}

\newtheorem{lem}{Lemma}
\newtheorem{theo}{Theorem}
\newtheorem{sled}{Corollary}
\newcommand{\proof}{\noindent\mbox{\bf Proof.~~}}

\begin{document}
\sloppy

\begin{center}
{\bf \large TRIGONOMETRIC POLYNOMIALS\\[3pt] DEVIATING THE LEAST FROM ZERO IN MEASURE\\[5pt]
AND RELATED PROBLEMS}

\ \

{\bf V.\,V.\,Arestov$^{1,2,}$\footnote{E-mail address:
Vitalii.Arestov@usu.ru}, A.\,S.\,Mendelev$^{1}$}

\ \

{\small {\it $^{1}$Ural State University, Yekaterinburg, Russia

$^{2}$Institute of Mathematics and Mechanics, Ural Division,
Russian Academy of Sciences,\\ Yekaterinburg, Russia}}
\end{center}

\begin{abstract}
We give a solution of the problem on trigonometric
polynomials~$f_n$ with the given leading harmonic $y\cos nt$
that deviate the least from zero in measure, more precisely,
with respect to the functional $\mu(f_{n})={\rm
mes}\{t\in[0,2\pi]:\ |f_{n}(t)|\geq 1\}$. For trigonometric
polynomials with a fixed  leading harmonic, we consider the
least uniform deviation from zero on a compact set and find the
minimal value of the deviation over compact subsets of the
torus that have a given measure. We give a solution of a
similar problem on the unit circle for algebraic polynomials
with zeros on the circle.
\end{abstract}


\section{Statement of the problem and preliminaries}

{\bf \thesection.1. Introduction.~ } Let $\mathcal{F}_{n}$ be the set of trigonometric
polynomials
\begin{equation}
\label{trig_polin} f_{n}(t)=\frac{a_{0}}{2}+\sum_{k=1}^{n}(a_{k}\cos kt+b_{k}\sin kt)
\end{equation}
of order $n\ge 0$ with real coefficients; in this paper, depending on the situation, we
consider these functions on the whole real line $\mathbb{R}$, on the period, i.e., a
segment of length~$2\pi$, or on the torus $\mathbb{T}$ which can be interpreted as a
segment of length~$2\pi$ (for example, the segment $[0,2\pi]$) with identified
end-points. On the set $\mathcal{F}_n$, we consider the functional
\begin{equation} \label{mera}
\mu(f_{n})={\mathrm{mes}\,}\{t\in\mathbb{T}:\ |f_{n}(t)|\geq 1\}
\end{equation}
whose value is the Lebesgue measure of the set of points of the torus at which an
absolute value of the polynomial $f_{n}\in \mathcal{F}_n$ is greater than or equal
to~$1$. For a fixed $y\ge 0$, we introduce the value
\begin{equation}
\label{basetask} \sigma_{n}(y)=\inf\{\mu\left(y\cos nt-f_{n-1}(t)\right):\ f_{n-1}\in
\mathcal{F}_{n-1}\}
\end{equation}
which can be interpreted as the value of the best approximation of the function $y\cos
nt$ by the set $\mathcal{F}_{n-1}$ of trigonometric polynomials of order $n-1$ with
respect to functional~\eqref{mera}. Value~\eqref{basetask} can be written in another
form. Let $\mathcal{F}_{n}(y)$ be the set of trigonometric polynomials of order~$n$ of
the form
$$
f_n(t)=y\cos nt+f_{n-1}(t),\qquad f_{n-1}\in \mathcal{F}_{n-1}.
$$
Then,
\begin{equation}
\label{basetproblem} \sigma_{n}(y)=\inf\{\mu\left(f_{n}\right):\ f_{n}\in
\mathcal{F}_{n}(y)\};
\end{equation}
this is a variant of the problem on polynomials that deviate the least from zero. It is
easily seen that problem~\eqref{basetask}--\eqref{basetproblem} is nontrivial only
for~\mbox{$y>1$}. The following assertion is valid for problem~\eqref{basetproblem}; in
this assertion and throughout the paper, we denote by~$T_n$ the Chebyshev polynomial of
the first kind which is specified by the formula $T_{n}(x)=\cos(n\arccos x)$ for $x\in
[-1,1]$.

\begin{theo}\label{maint}
For $y>1$ and $n\ge 1$, the following equality is valid:
\begin{equation}\label{Eq}
\sigma_{n}(y)=4\arccos\frac{1}{y^{\frac{1}{2n}}}.
\end{equation}
Moreover, for $k=0,1,\ldots 2n-1$, the polynomial
\begin{equation}
\label{extr-polinom-BAG} f_{n}(t)=f_{n,k}(t)=(-1)^kT_{n}\left( y^{\frac{1}{n}}\cos
\left(t-\frac{\pi k}{n}\right)-y^{\frac{1}{n}}+1\right)
\end{equation}
belongs to the set $\mathcal{F}_n(y)$ and it is an extremal polynomial in
problem~\eqref{basetproblem}~$($i.e., the polynomial that deviates the least from zero$)$
and only such polynomials solve problem~\eqref{basetproblem}.
\end{theo}

A.\,S.\,Mendelev announced this result in 2000 in abstracts of his talk~\cite{m}. The
proof of Theorem~\ref{maint} is published in the present paper for the first time. In
1998, A.\,S.\,Mendelev and M.\,S.\,Plotnikov~\cite{b1} proved assertion~\eqref{Eq} for
large values of~$y$; more precisely, for
$$
y\geq \frac{1}{\sin^{2n}\dfrac{\pi}{2(2n+1)}}.
$$

In 1992, A.\,G.\,Babenko~\cite{BAG} studied the least constant $\beta_{n}$ in the
inequality
\begin{equation}\label{babenko-ineq}
\mu(a_n\cos nt+b_n\sin nt)\le \beta_{n}\, {\mu(f_n) },\qquad f_n\in \mathcal{F}_{n},
\end{equation}
on the set $\mathcal{F}_n$ of trigonometric polynomials~\eqref{trig_polin}; he obtained
the following estimates for $\beta_n$:
\begin{equation}
\label{babenko_i} \sqrt{2n} \leq \beta_{n} \leq n\sqrt{2},\qquad n \geq 1.
\end{equation}
Below (see Theorem~\ref{Ineq-Babenko}), as a consequence of Theorem~\ref{maint}, we find
the value $\beta_{n}$; namely, we show that $\beta_{n}=\sqrt{2n}$. Thus, it turned out
that the lower bound in~\eqref{babenko_i} is true.

In this paper, in connection with the investigation of problem~\eqref{basetproblem}, we
discuss several other related extremal problems for trigonometric polynomials on the
torus~$\mathbb{T}$ and for algebraic polynomials on the unit circle $\Gamma$ of the
complex plane. In particular, we produce the following results.

1) For algebraic polynomials with zeros on the unit circle of
the complex plane and with the unit leading coefficient, we
consider the least uniform deviation from zero on a compact set
and find the minimal value of the deviation over compact
subsets of the circle that have a given measure.

2) For trigonometric polynomials with a fixed  leading
harmonic, we consider the least uniform deviation from zero on
a compact set and find the minimal value of the deviation over
compact subsets of the torus that have a given measure.

The main part of the results of this paper were stated without proofs
in~\cite{Arestov-Mendelev-DAN}.

\ \

{\bf \thesection.2. A restriction of the class of polynomials.~} We need certain known
facts about trigonometric polynomials~\eqref{trig_polin}; for further actions, it is
sufficient to consider only polynomials whose order is equal to~$n$, i.e., such that
$a_n^2+b_n^2>0$. For a trigonometric polynomial~$f_n$ of order $n\ge 1$ with real
coefficients, the following formula is valid:
\begin{equation}
\label{a1} f_n(t)=e^{-int} P_{2n}(e^{it}),\qquad t\in \mathbb{R},
\end{equation}
where
\begin{equation}
\label{a2} P_{2n}(z)=\sum_{\nu=0}^{2n}u_\nu z^\nu
\end{equation}
is an algebraic polynomial of degree~$2n$ whose coefficients have the properties
\begin{equation} \label{a3}
u_{2n-\nu}={\overline u}_\nu,\qquad 0\le \nu\le 2n.
\end{equation}
Conversely, if the coefficients of polynomial~\eqref{a2} satisfy conditions~\eqref{a3},
then formula~\eqref{a1} specifies a trigonometric polynomial of order~$n$ with real
coefficients. In this case, in particular,
\begin{equation} \label{a4}
u_{2n}=\overline{u}_0=\frac{a_n-ib_n}{2}.
\end{equation}

Condition~\eqref{a3} means that the following formula is valid for polynomial~\eqref{a2}:
\begin{equation}
\label{a5} z^{2n}\overline{P_{2n}\left({\overline z}\,^{-1}\right)}=P_{2n}(z),\qquad
z\in\mathbb{C},\qquad z\ne 0.
\end{equation}
Hence, polynomial~\eqref{a2} can  be written in the form
\begin{equation}
\label{complex}
P_{2n}(z)=\frac{a_{n}-ib_{n}}{2}\left[\prod_{k=1}^{{l}}(z-z_{k})\left(z-\overline{z_{k}}\
^{-1}\right)\right] \prod_{j={2{l}}+1}^{2n}\left(z-e^{i\phi_{j}}\right).
\end{equation}
In this representation, the first product corresponds to ${2{l}}$ complex zeros of the
polynomial~$f_n$ and $0<|z_k|<1,$~ $1\le k\le {l}$; if~$f_n$ has no complex zeros
$({l}=0)$, then this representation is absent. The second product in~\eqref{complex}
corresponds to $2(n-{l})$ real zeros of the polynomial~$f_n$; if all zeros of the
polynomial~$f_n$ are complex, then the second product in~\eqref{complex} is absent.

However, there exist polynomials of the form~\eqref{complex} without property~\eqref{a5}.
Substituting expression~\eqref{complex} into~\eqref{a5}, we ascertain that~\eqref{a5}
holds only in the case if
\begin{equation} \label{a6}
2\theta_n+2\sum_{k=1}^{{l}}\varphi_k+\sum_{j={2{l}}+1}^{2n}\phi_j=2\pi N,\qquad N\in
\mathbb{Z},
\end{equation}
where $\theta_n$ is an argument of the coefficient $a_n-ib_n$, and $\varphi_k$ are
arguments of the zeros~$z_k$, ~$1\le k\le {{l}}$, of the polynomial $P_{2n}$. Thus,
relation~\eqref{a6} is a necessary and sufficient condition for
polynomial~\eqref{complex} to have property~\eqref{a3}, and so to generate, by
formula~\eqref{a1}, a trigonometric polynomial of order~$n$ with real coefficients; in
addition, the leading harmonic of the polynomial has the form $a_n\cos nt+b_n\sin nt$. A
more detailed information on the facts presented here can be found, for example,
in~\cite[Sect. VI,\ Subsect. 2]{P-S}.

Let us discuss the representation of polynomials $f_n\in \mathcal{F}_n(y)$ (for $y>0$ and
$n\ge 1$) in more details. In this case, $a_n=y>0$ and $b_n=0$; therefore, $\theta_n=0$.
Polynomial~\eqref{complex} and condition~\eqref{a6} take the form
\begin{equation}
\label{complex-y}
P_{2n}(z)=\frac{y}{2}\left[\prod_{k=1}^{{l}}(z-z_{k})\left(z-\overline{z_{k}}\
^{-1}\right)\right] \prod_{j={2{l}}+1}^{2n}\left(z-e^{i\phi_{j}}\right),
\end{equation}
\begin{equation} \label{a6-y} 2\sum_{k=1}^{{l}}\varphi_k+\sum_{j={2{l}}+1}^{2n}\phi_j=2\pi N,\qquad N\in \mathbb{Z}.
\end{equation}
By formulas~\eqref{a1} and~\eqref{complex-y}, the following equality is valid:
\begin{equation}
\label{real_T} f_{n}(t)=e^{-int}\frac{y}{2}\prod_{k=1}^{{l}}\left(\left(e^{it}-{r_k}
e^{i{\varphi_k}}\right)\left(e^{it}-\frac{e^{i{\varphi_k}}}{{r_k}}\right)\right)
\prod_{j={2{l}}+1}^{2n}(e^{it}-e^{i\phi_{j}}),
\end{equation}
where ${r_k}=|z_k|$ and $\varphi_k=\arg z_k$ are the modulus and argument of the
zero~$z_k$, ~$1\le k\le {{l}}$, respectively. Let us simplify the right-hand side of
representation~\eqref{real_T}. For the multipliers from the second product, we have
$$
e^{it}-e^{i\phi_{j}}=ie^{i\frac{t+\phi_j}{2}}\left(2\sin\frac{t-\phi_{j}}{2}\right).
$$
Let us transform the multipliers from the first product as follows:
$$
\left(e^{it}-{r_k}
e^{i{\varphi_k}}\right)\left(e^{it}-\frac{e^{i{\varphi_k}}}{{r_k}}\right)=
e^{2it}-\left({r_k}+\frac{1}{{r_k}}\right)e^{i(t+{\varphi_k})}+e^{2i{\varphi_k}}=
$$
$$
=-2e^{i(t+\varphi_k)}(A_k-\cos(t-{\varphi_k})),\qquad
A_k=\dfrac{1}{2}\left({r_k}+\dfrac{1}{{r_k}}\right)>1.
$$
Thus, the following representation is valid for a polynomial $f_n\in \mathcal{F}_n(y)$:
\begin{equation}
\label{T-predstavlenie}
f_{n}(t)=(-1)^{n+N}\,\frac{y}{2}\,\prod_{k=1}^{{l}}\left(2(A_{k}-\cos(t-{\varphi_k}))\right)
\prod_{j={2{l}}+1}^{2n}\left(2\sin\frac{t-\phi_{j}}{2}\right),
\end{equation}
where $A_{k}>1$ for $1\le k\le {{l}}$; $\varphi_k\in \mathbb{R}$ for $1\le k\le l$; and
$\phi_j\in \mathbb{R}$ for ${2{l}}+1\le j\le 2n$. We recall that, in addition,
condition~\eqref{a6-y} is valid.

In what follows, the set $\mathcal{F}_n^{\,real}(y)$ of trigonometric polynomials from
$\mathcal{F}_n(y)$ all zeros of which are real will play an important role. By~\eqref{a1}
and~\eqref{complex-y}, for a polynomial $f_n\in \mathcal{F}_n^{\,real}(y)$, we have
\begin{equation}
\label{T-real} f_n(t)=e^{-int}\,\frac{y}{2} \,P_{2n}(e^{it}),\qquad
P_{2n}(z)=\prod_{j=1}^{2n}\left(z-e^{i\phi_{j}}\right);
\end{equation}
here, $\phi_{j}\in \mathbb{R}$, ~$1\le j\le 2n$, and the following condition holds:
\begin{equation}
\label{a6-y-refl} \sum_{j=1}^{2n}\phi_j=2\pi N,\qquad N\in \mathbb{Z}.
\end{equation}
Formula~\eqref{T-real} also implies the representation
\begin{equation}
\label{T-real-predstavlenie} f_{n}(t)=(-1)^{n+N}\,\frac{y}{2}\,g_{2n}(t),\qquad
g_{2n}(t)=\prod_{j=1}^{2n}\left(2\sin\frac{t-\phi_{j}}{2}\right);
\end{equation}
here, $\phi_j$, ~$1\le j\le 2n$, are (real) zeros of the polynomial~$f_n$. From the above
reasonings, it is clear that condition~\eqref{a6-y-refl} is necessary and sufficient for
function~\eqref{T-real} or, that is the same,~\eqref{T-real-predstavlenie} to be a
polynomial from~$\mathcal{F}_n^{\,real}(y)$.

Justifying results of~\cite{b1}, A.\,S.\,Mendelev and M.\,S.\,Plotnikov obtained the
following assertion.

\begin{lem}\label{lem1}
For any $ n\ge 1$ and $y>1$, an extremal polynomial exists in
problem~\eqref{basetproblem}; all zeros of the extremal polynomial are real.
\end{lem}

\proof Assume that a polynomial $f_n\in \mathcal{F}_n(y)$ has at least one complex zero;
i.e., representation~\eqref{T-predstavlenie} contains at least one multiplier
$A_{k}-\cos(t-\varphi_{k})$. Since $A_k>1$, the following inequality is valid for
all~$t$:
$$
A_k-\cos(t-\varphi_{k})>1-\cos(t-\varphi_{k})= 2\sin^{2}\frac{t-\varphi_{k}}{2}.
$$
Let us consider the function
\begin{equation}
\label{a7} \widehat{f}_n(t)=e^{-int}\,\frac{y}{2}
\,\prod_{k=1}^{{l}}\left(e^{it}-e^{i\varphi_k}\right)^2\prod_{j=
{2{l}}+1}^{2n}\left(e^{it}-e^{i\phi_{j}}\right),
\end{equation}
which is a trigonometric polynomial of order~$n$. The polynomial~$f_n$ satisfies
condition~\eqref{a6-y-refl}; hence, function~\eqref{a7} is a trigonometric polynomial
with real zeros; more precisely, $\widehat{f}_n\in \mathcal{F}_n^{\,real}(y)$. Absolute
values of the polynomials~$f_n$ and $\widehat{f}_n$ are connected by the inequality
$|\widehat{f}_n(t)|\le|{f}_n(t)|,\ t\in \mathbb{R}$; moreover, the strict inequality
$|\widehat{f}_n(t)|<|{f}_n(t)|$ holds at points $t\in \mathbb{R}$ distinct from real
zeros $\{\phi_j\}_{j={{2{l}}+1}}^{2n}$ of the polynomial~$f_n$. Hence,
$\mu(\widehat{f}_n)<\mu(f_n)$.

At this stage, we, in particular, have proved that, in~\eqref{basetproblem}, it is
necessary to restrict our attention to polynomials $f_{n}\in
\mathcal{F}_{n}^{\,real}(y)$; consequently, the following equality holds:
\begin{equation}
\label{realproblem} \sigma_{n}(y)=\inf\{\mu\left(f_{n}\right):\ f_{n}\in
\mathcal{F}_{n}^{\,real}(y)\}.
\end{equation}
To complete the proof of the lemma, it remains to show that the infimum in the right-hand
side of~\eqref{realproblem} is reached. For polynomials $f_{n}\in
\mathcal{F}_{n}^{\,real}(y)$, the value $\mu(f_n)$ is a function of $2n$ (real) zeros of
the polynomial $ f_{n}\in \mathcal{F}_{n}^{\,real}(y)$. We will use the same symbol~$\mu$
to denote this function; thus, $\mu (f_n)=\mu(\phi_1,\dots,\phi_{2n})$. The zeros
$\{\phi_j\}_{j=1}^{2n}$ are related by condition~\eqref{a6-y-refl}. We can assume that
$0\le \phi_j\le 2\pi,\ 1\le j\le 2n-1$, and $\phi_{2n}=-\sum_{j=1}^{2n-1}\phi_j$. This
set of points $\phi=\{\phi_j\}_{j=1}^{2n}$ will be denoted by $\Pi_{2n}$.
Relation~\eqref{realproblem} can be rewritten in the form
\begin{equation}
\label{Mproblem} \sigma_{n}(y)=\inf\{\mu(\phi_1,\dots,\phi_{2n}):\
\phi=\{\phi_j\}_{j=1}^{2n}\in \Pi_{2n}\}.
\end{equation}
It is easily seen that the function $\mu$ continuously depends on the point
$\phi=\{\phi_j\}$; in addition, $\Pi_{2n}$ is a compact subset of the
space~$\mathbb{R}^{2n}$. Therefore, the infimum in~\eqref{Mproblem} or, that is the same,
in~\eqref{realproblem} is reached. The proof of the lemma is completed.

\ \

{\bf \thesection.3. A restatement and expansion of the initial problem.} Lemma~\ref{lem1}
reduces initial problem~\eqref{basetproblem} to a more clear problem of minimization of a
(continuous) function of several real variables. Let~$m$ be natural; in the sequel,
studying problem~\eqref{basetproblem}, we take~$m=2n$. We denote by
$\mathfrak{P}_m(\Gamma)$ the set of algebraic polynomials
\begin{equation}
\label{Pm} P_m(z)=\prod_{j=1}^m(z-e^{i\phi_j}),\qquad z\in \mathbb{C},
\end{equation}
of order~$m$ with the unit leading coefficient, all~$m$ zeros of which belong to the unit
circle $\Gamma=\{e^{it}:\ t\in[0,2\pi]\}$. Every such polynomial is uniquely defined by
the point $\phi=(\phi_1,\dots,\phi_m)\in{{\mathbb R}^m}$. On the unit circle,
polynomial~\eqref{Pm} is representable in the form
\begin{equation}
\label{G-s} P_m(e^{it})=e^{i\frac{m}{2}t}e^{i\frac{\Phi}{2}}i^mg_m(t),\qquad
\Phi=\sum_{j=1}^m\phi_j,
\end{equation}
where
\begin{equation}
\label{bb2} g_m(t)=g_m(t;\phi)=
\prod_{j=1}^{m}\left(2\sin\frac{t-\phi_{j}}{2}\right),\qquad t\in {\mathbb R}.
\end{equation}
The set of functions~\eqref{bb2} will be denoted by $\mathcal{G}_m$.

Let us consider the quantity
\begin{equation}
\label{c-norm-S} {h(}m)=\min\{ \|g\|_{C_{2\pi}}:\ g\in \mathcal{G}_m\}=\min\{
\|P_m\|_{C(\Gamma)}:\ P_m\in \mathfrak{P}_m(\Gamma)\}
\end{equation}
of the least value of the uniform norm of polynomials~\eqref{bb2} on the real line or,
that is the same, of the uniform norm of polynomials~\eqref{G-s} on the unit circle. It
is well known that
\begin{equation}
\label{c-norm-S-v} h(m)=2, \qquad m\ge 1.
\end{equation}
Besides, it  can be easily verified. Indeed, polynomial~\eqref{Pm} has the form
\begin{equation}
\label{Pm-coef} P_m(z)=\sum_{k=0}^mc_kz^k;
\end{equation}
here, $c_m=1$, and $c_0=e^{i\Psi}$, ~$\Psi=m\pi+\sum_{j=1}^m\phi_j $. For a fixed $\Psi$,
on the set of algebraic polynomials $F_m\in\mathcal{P}_m$ of order $m\ge 1$, let us
consider the linear functional
$$
\Sigma_m(F_m)=\frac{1}{m}\sum_{{l}=0}^{m-1}F_m(e^{i\theta_{l}}),\qquad{\mbox{where}}\qquad
\theta_{l}=\frac{\Psi+2\pi{l}}{m} ,\qquad 0\le{l}\le m-1.
$$
For the polynomials $p\,_k(z)=z^k$, we have
$$
\Sigma_m(p\, _k)=\frac{1}{m}\sum_{{l}=0}^{m-1}e^{ik\theta_{l}}=e^{i\frac{k\Psi}{m}}
\frac{1}{m}\sum_{{l}=0}^{m-1}e^{i\frac{2\pi k{l}}{m}}.
$$
Hence, we see that $\Sigma_m(p_k)=0,\, 1\le k\le {m-1}$; $\Sigma_m(1)=1$;
$\Sigma_m(p_m)=e^{i\Psi}$. Therefore, $ \Sigma_m(P_m)=2e^{i\Psi}$ for
polynomial~\eqref{Pm}. On the other hand, the estimate $|\Sigma_m(P_m)|\le
\|P_m\|_{C(\Gamma)}$ is valid. Consequently, $h(m)\ge 2$. The polynomial $P_m(z)=z^m+1$
provides the inverse estimate. Thus, assertion~\eqref{c-norm-S-v} really holds.

For a parameter $h,\ 0\le h\le 2=h(m)$, we set
\begin{equation}
\label{task2-m-P} \delta_{m}({h})=\inf\{{\mathrm{mes}\,}\{t\in\mathbb{T}:\
|P_m(e^{it})|\geq h\}:\ P_m\in\mathfrak{P}_m(\Gamma)\}.
\end{equation}
Relations~\eqref{G-s} and~\eqref{bb2} imply also that
\begin{equation}
\label{task2-m-G} \delta_{m}({h})=\inf\{{\mathrm{mes}\,}\{t\in\mathbb{T}:\ |g_m(t)|\geq
h\}:\ g_m\in\mathcal{G}_m\}=
\end{equation}
\begin{equation}
\label{task2-m}=\inf\{{\chi}_{m}(\phi\, ; h):\ \phi=(\phi_1,\dots,\phi_{m})\in {{\mathbb
R}^m}\},
\end{equation}
where
\begin{equation}
\label{mera-m} {\chi}_{m}(\phi;h)={\chi}_{m}(g_m;h)={\mathrm{mes}\,}\{t\in[0,2\pi]:\
|g_m(t;\phi)|\geq h\}
\end{equation}
is a function of the point $\phi=(\phi_1,\dots,\phi_{m})\in {{\mathbb R}^m}$.
Representation~\eqref{task2-m} means that (for a fixed $h\in[0,2]$) the value
$\delta_m(h)$ can be interpreted as the minimum of a (continuous) function
${\chi}_{m}(\phi)={\chi}_{m}(\phi;h)$ of~$m$ variables. A considerable part of this paper
is devoted to studying the value~$\delta_{m}({h})$. In the sequel, depending on the
situation, it will be convenient for us to use one of the three representation
forms~\eqref{task2-m-P}--\eqref{task2-m} for the value~$\delta_m(h)$.

In the case $m=1$, function~\eqref{bb2} takes the form
\begin{equation}
\label{bb2-1} g_1(t)=2\sin\frac{t-\phi_{1}}{2}.
\end{equation}
For any such function, ${\mathrm{mes}\,}\{t\in\mathbb{T}:\ |g_1(t)|\geq h\}=4\arccos
(h/2)$. Consequently, the following formula holds for $m=1$:
\begin{equation}\label{delta_h_m1}
\delta_1({h})=4\arccos\left(\frac{h}{2}\right),\qquad 0\le h\le 2.
\end{equation}
Moreover, any polynomial~\eqref{bb2-1} is extremal in~\eqref{task2-m-G}; so, any
polynomial $P_1(z)=z-\zeta$, whose zero $\zeta$ satisfies the condition $|\zeta|=1$, is
extremal in~\eqref{task2-m-P}.

\begin{lem}\label{lem_delta}
For $m\ge 1$, the following assertions are valid$:$

$1)$ for any $h,\ 0\le h\le 2$, there exists an extremal point $\phi=\phi^*\in
\mathbb{R}^m$, at which an infimum in~\eqref{task2-m} is reached$;$

$2)$ for the extreme values of~$h$, we have $\delta_{m}(0)=2\pi$, ~$\delta_{m}(2)=0$;

$3)$ the value $\delta_{m}({h})$~$($strictly$)$ decreases with respect to~$h\in [0,2]$.
\end{lem}

\proof It is sufficient to consider~$m\ge 2$. If $h=0$, then, for any
function~\eqref{bb2}, value~\eqref{mera-m} is equal to~$2\pi$; so $\delta_{m}(0)=2\pi$.
Let us discuss the case~$h=2$. For example, the polynomial $ P_m(z)=z^m+1 $ belongs to
the set $\mathfrak{P}_m(\Gamma)$; for this polynomial, the set $\{t\in\mathbb{T}:\
|P_m(e^{it})|\geq 2\}$ consists of~$m$ points; thus, its measure is zero.
Therefore,~$\delta_{m}(2)=0$.

The existence of an extremal point for $0<h<2$ in~\eqref{task2-m} (and so, of extremal
functions in~\eqref{task2-m-P} and~\eqref{task2-m-G}) can be easily justified with the
help of the arguments used in the proof of Lemma~\ref{lem1}.

Finally, let us prove the monotonicity of the value $\delta_m(h)$ with respect
to~$h\in[0,2]$. Let $0\le h_1<h_2\le 2$. We denote by $g_m^{(1)}$ the polynomial from
$\mathcal{G}_m$ on which the infimum in~\eqref{task2-m-G} is reached for~$h=h_1$. The
strict inequality $ {\mathrm{mes}\,}\{t\in\mathbb{T}:\ |g_m^{(1)}(t)|\geq
h_2\}<{\mathrm{mes}\,}\{t\in\mathbb{T}:\ |g_m^{(1)}(t)|\geq h_1\}$ holds. This implies
that $\delta_{m}({h_2})<\delta_{m}({h_1}),\ 0\le h_1<h_2\le 2$. The lemma~is proved.

We denote by $\mathbb{H}=\mathbb{H}_{m}$ the hyperplane of points
$\phi=(\phi_1,\dots,\phi_{m})\in {\mathbb R}^{m}$ satisfying the condition
\begin{equation} \label{beta-m} \sum_{j=1}^{m}\phi_j=0.
\end{equation}
This hyperplane is orthogonal to the vector $\mathcal{E}=\mathcal{E}_m=(1,1,\ldots,1) \in
{{\mathbb R}^m}$. Let us ascertain that, in~\eqref{task2-m}, we can restrict our
attention to points $\phi=(\phi_1,\dots,\phi_{m})\in\mathbb{H}_m$; more precisely, that
the following formula holds:
\begin{equation}
\label{task2-m-beta} \delta_{m}({h})=\inf\{{\chi}_{m}(\phi\, ; h):\
\phi=(\phi_1,\dots,\phi_{m})\in \mathbb{H}_m\},\qquad 0\le h\le 2.
\end{equation}
Indeed, let $\phi=(\phi_1,\dots,\phi_{m})\in {{\mathbb R}^m}$. We set
$\phi_0=(\phi_1+\dots+\phi_m)/m$. Then, the point $
\overline\phi=\phi-\phi_0\mathcal{E}=(\phi_1-\phi_0,\dots,\phi_{m}-\phi_0) $ belongs to
the hyperplane $\mathbb{H}_m$. It is easily seen that the equality
$\chi_m(\phi,h)=\chi_m(\overline\phi,h)$ holds. Hence, relation~\eqref{task2-m-beta}
follows.

By Lemma~\ref{lem1} and formula~\eqref{T-real-predstavlenie}, the following assertion is
valid.

\begin{sled} \label{sigma-delta}
For $n\ge 1$ and $y>1$, values~\eqref{basetproblem} and~\eqref{task2-m-P} are related as
follows:
\begin{equation}
\label{problem-2n-var2} \sigma_{n}(y)=\delta_{2n}({h}),\qquad {h}=\frac{2}{y}.
\end{equation}
\end{sled}

\section{A problem equivalent to problem~\eqref{task2-m}\\ and its investigation}
\setcounter{equation}{0}

{\bf \thesection.1. An equivalent problem.~ } In this section, we will study a problem
equivalent to problem~\eqref{task2-m}. For natural~$m$ and real $h\ge 0$, we introduce
the set
\begin{equation}
\label{b1} {{\mathbb V}}={\mathbb V}(h)=\left\{x=(x_1,x_2,\ldots,x_m)\in {{\mathbb R}^m}
:\ \left|\prod\limits_{j=1}^m s(x_j)\right|\geq h \right\}\subset {{\mathbb R}^m};
\end{equation}
here and subsequently,
$$s(t)=2\sin \frac{t}{2},\qquad t\in\mathbb{R}.$$
Set~\eqref{b1} is nonempty if and only if~$h\le 2^m$. However, as will be seen below, we
are interested only in values $0<h<2$. To a number $a\in \mathbb{R}$ and point
$\phi=(\phi_1,\dots,\phi_m)\in {{\mathbb R}^m}$ we assign the set
\begin{equation}
\label{b2} \Upsilon_m(\phi,h;a)=[\phi+a\mathcal{E}, \phi+(a+2\pi)\mathcal{E}] \cap
{{\mathbb V}(h)}\subset {{\mathbb R}^m},
\end{equation}
which is the intersection of the segment $[\phi+a\mathcal{E}, \phi+
(a+2\pi)\mathcal{E}]\subset{{\mathbb R}^m}$ with set~\eqref{b1}; here,
$\mathcal{E}=(1,1,\ldots,1) \in {{\mathbb R}^m}$. We are interested in the linear measure
${\mathrm{mes}}_1(\Upsilon_m(\phi,h;a))$ of this set. The subsequent considerations will
show that this measure is independent of the parameter~$a$. We set
\begin{equation}
\label{b3} \Delta_{m}(h)=\inf\{ {\mathrm{mes}}_1\left([\phi+a\mathcal{E} ,\phi
+(a+2\pi)\mathcal{E}] \cap {{\mathbb V}}(h)\right):\ \phi \in {{\mathbb R}^m} \}.
\end{equation}

As will be shown below (see Lemma~\ref{lemDelta}), problems~\eqref{b3}
and~\eqref{task2-m} are closely interrelated; namely, the following equaity holds:
\begin{equation} \label{b4}
\Delta_{m}(h)=\sqrt{m}\cdot \delta_{m}(h),\qquad 0<h<2.
\end{equation}
To make sure in this, let us compare the linear measure
${\mathrm{mes}}_1(\Upsilon_m(\phi,h;a))$ of set~\eqref{b2} and the measure
${\mathrm{mes}}(\upsilon_m(\phi,h;a))$ of the set
\begin{equation}
\label{b5} \upsilon_m(\phi,h;a)=\{t \in [a,a+2\pi]: ~|g_m(t;\phi)| \geq h\}\subset
{{\mathbb R}}
\end{equation}
constructed by the function
\begin{equation}
\label{b6} g_m(t)=g_m(t;\phi)=\prod_{j=1}^{m}\left(2\sin\frac{t+\phi_{j}}{2}\right)=
\prod_{j=1}^{m}s({t+\phi_{j}}).
\end{equation}
Note that, here, in comparison with the previous section, the zeros~$\phi_j$ have the
reversed signs.

\begin{lem}\label{lem2}
For a point $\phi=(\phi_1,\dots,\phi_m) \in {{\mathbb R}^m}$ and parameter $h,\ 0<h<2$,
the following assertions are valid$:$ $1)$ there holds the equality
\begin{equation} \label{b7}
{\mathrm{mes}}_1([\phi+a\mathcal{E}, \phi+(a+2\pi)\mathcal{E}] \cap {{\mathbb
V}(h)})=\sqrt{m}\cdot {\mathrm{mes}\,}\{t \in [a,a+2\pi]: ~|g_m(t;\phi)| \geq h\};
\end{equation}
$2)$ both sets~\eqref{b5} and~\eqref{b2} consist of the same number of
segments~$($probably, degenerating to a point$)$ whose lengths are directly proportional
with the coefficient~$\sqrt{m};$ $3)$ the measures of sets~\eqref{b5} and~\eqref{b2} are
independent of the parameter~$a\in \mathbb{R}$.
\end{lem}

\proof Let us consider the linear vector-function $\alpha(t)=\phi+t\mathcal{E}=
(\alpha_1(t), \ldots, \alpha_m(t))$, where $\alpha_j(t)=\phi_j+t$, $j=1,\ldots,m$,
~$t\in[a,a+2\pi]$. This function is a bijection of the segment $[a,a+2\pi]$ onto the
segment $[\phi+a\mathcal{E} ,\phi+(a+2\pi)\mathcal{E}]$. An interval $X\subset
[a,a+2\pi]$ is mapped by the function $\alpha$ onto an interval
$\mathcal{X}=\alpha(X)\subset [\phi+a\mathcal{E} ,\phi+(a+2\pi)\mathcal{E}]$ of the same
type; moreover, it is easily seen that the measure ${\mathrm{mes}\,}(X)$ of an interval
$X\subset [a,a+2\pi]$ and the linear measure ${\mathrm{mes}}_1(\mathcal{X})$ of the
interval $\mathcal{X}=\alpha(X)$ are related by the equality
${\mathrm{mes}}_1(\mathcal{X})=\sqrt{m}\cdot {\mathrm{mes}\,}(X)$. It easily follows
that, for any measurable subset $X\subset [a,a+2\pi]$, its image
$\mathcal{X}=\alpha(X)\subset [\phi, \phi+2\pi\mathcal{E}]$ is also measurable and the
(linear) measures of these sets are related by the equality
${\mathrm{mes}}_1(\mathcal{X})=\sqrt{m}\cdot {\mathrm{mes}\,}(X)$.

Let us ascertain that set~\eqref{b2} is the image of set~\eqref{b5} under the
mapping~$\alpha$; i.e.,
\begin{equation}
\label{b8} \Upsilon_m(\phi,h;a)=\alpha(\upsilon_m(\phi,h;a)).
\end{equation}
The fact that the point $\alpha(t)=\phi+t\mathcal{E}$ belongs to set~\eqref{b2} means
that this point lies on the segment $[\phi+a\mathcal{E}, \phi+ (a+2\pi)\mathcal{E}]$ and
in the set ${{\mathbb V}(h)}$ simultaneously. The first fact means that $t\in
[a,a+2\pi]$. By definition~\eqref{b1}, the fact that the point $\alpha(t)$ belongs to the
set ${{\mathbb V}(h)}$ means that $\prod_{k=1}^m |s(\alpha_k(t))|\geq h$. However,
$$
\prod\limits_{j=1}^m s(\alpha_j(t))=\prod\limits_{j=1}^m s(\phi_j+t)=
\prod\limits_{j=1}^{m}\left(2\sin\frac{t+\phi_j}{2}\right)=g_m(t;\phi).
$$
Thus, $\alpha(t)\in \Upsilon_m(\phi,h;a)$ if and only if $t\in \upsilon_m(\phi,h;a)$.
Assertion~\eqref{b8} is proved.

Set~\eqref{b5} consists of a finite number of segments (some of which can degenerate to a
point). Set~\eqref{b2} has the same structure. Both the sets are measurable and their
measures are related by equality~\eqref{b7}.

The measure of set~\eqref{b5} as well as, by equality~\eqref{b7}, the measure of
set~\eqref{b2} are independent of the number~$a$. The lemma is proved.

\begin{lem}\label{lemDelta}
For $m\ge 1$ and $ 0<h<2$, the following assertions are valid for
problems~\eqref{task2-m} and~\eqref{b3}$:$ $1)$~equality~\eqref{b4} holds$;$ $2)$~there
exists an extremal point $\phi=\phi^* \in \mathbb{R}$ at which infimums
in~\eqref{task2-m} and~\eqref{b3} are reached$;$ this point has property~$2$ from
Lemma~$\ref{lem2};$ $3)$~value~\eqref{b3} is independent of the parameter~$a\in
\mathbb{R}$.
\end{lem}

\proof Equality~\eqref{b4} is a consequence of the previous lemma. By Lemma~\ref{lem2},
it is also sufficient to justify the existence of an extremal point in~\eqref{task2-m};
this have been done in Lemma~\ref{lem_delta}. Lemma~\ref{lemDelta} is proved.

Our immediate aim is to  ascertain that an extremal set in problem~\eqref{b3} (i.e.,
set~\eqref{b2} for the extremal point $\phi=\phi^*$ of problem~\eqref{b3}) consists of
one segment and $m-1$ points. This, Lemma~\ref{lem2} and Lemma~\ref{lemDelta} will imply
that the set $\{t\in\mathbb{T}:\,|g_{m}(t)|\geq h\}$ for the extremal polynomial~$g_m$
in~\eqref{task2-m-G} or, that is the same, the set
$\{t\in\mathbb{T}:\,|P_{m}(e^{it})|\geq h\}$ for the extremal polynomial~$P_m$
in~\eqref{task2-m-P} also consist of one segment and $m-1$ points.

\ \

{\bf \thesection.2. Properties of the set ${\mathbb V}$.~ } Starting with the set
${\mathbb V}={\mathbb V}(h)$ defined by~\eqref{b1}, we introduce the set ${\mathbb
V}_0={\mathbb V}_0(h)={{\mathbb V}}(h) \cap (0, 2\pi)^m$. The function~$s$ has the
property $s(2l \pi)=0$ for $l\in {\mathbb Z}$; therefore,
\begin{equation}\label{b10}
{\mathbb V}_0(h)={{\mathbb V}}(h) \cap (0, 2\pi)^m={{\mathbb V}}(h) \cap [0,
2\pi]^m,\qquad h>0.
\end{equation}
Evidently, for any $h,\ 0\le h\le 2^m$, the set ${\mathbb V}_0(h)$ is nonempty; and, for
$0\le h<2^m$, this set consists of more than one point. In addition, the sets $ {\mathbb
V}_0(h)$ decreases with respect to~$h$; more precisely,
\begin{equation}\label{b15}
{\mathbb V}_0(h_2) \subset {\mathbb V}_0(h_1),\qquad 0\le h_1\le h_2\le 2^m.
\end{equation}

Along with ${\mathbb V}_0$, we consider the sets ${\mathbb V}_k={\mathbb V}_k(h)={\mathbb
V}_0+2 \pi k=\{x+2\pi k:\ x\in {\mathbb V}_0\}$ that are shifts of set~\eqref{b10} by the
vectors $2 \pi k,\ k \in {{\mathbb Z}^m}$.

\begin{lem}\label{lemV}
For $m\ge 1$ and $0<h<2^m$, the following assertions are valid$:$
\begin{equation}\label{b11}
\text{the sets}~ \{{\mathbb V}_k:\ k \in {{\mathbb Z}^m}\}~ \text{are pairwise disjoint},
\end{equation}
\begin{equation}\label{b12}
{{\mathbb V}}=\bigcup\limits_{k \in {{\mathbb Z}^m}} {\mathbb V}_k,
\end{equation}
\begin{equation}\label{b13}
{\mathbb V}_0~ \text{is compact},
\end{equation}
\begin{equation}\label{b14}
{\mathbb V}_0~ \text{is strictly convex},
\end{equation}
\begin{equation}\label{b13a}
{\mathbb V}_0~ \text{has a nonempty interior; i.e., it is a body in}~ {\mathbb R}^m.
\end{equation}
\end{lem}

\proof Property~\eqref{b11} is evident.

The set $\mathbb{V}$ (see definition~\eqref{b1}) can be written in the form
$$
{{\mathbb V}}= \left\{x=(x_1,x_2,\ldots,x_m)\in {{\mathbb R}^m} :\ w(x)\geq h \right\},
$$
where
\begin{equation} \label{fun-w}
{w}(x)={w}(x_1,\dots,x_m)=\prod\limits_{j=1}^m |s(x_j)|,\qquad
x=(x_1,\dots,x_m)\in\mathbb{R}^m.
\end{equation}
The function $|s(t)|=2|\sin (t/2)|$ is continuous and $2\pi$-periodic on the real line.
Therefore, function~\eqref{fun-w} is continuous on $\mathbb{R}^m$ and $2\pi$-periodic
with respect to every variable; more precisely, $w(x+2\pi k)=w(x)$ for all
$x\in\mathbb{R}^m$ and $k\in\mathbb{Z}^m$. This, in particular, implies
property~\eqref{b12}.

Since function~\eqref{fun-w} is continuous on ${{\mathbb R}^m},$ the set ${{\mathbb V}}$
is closed. By representation~\eqref{b10}, the set ${\mathbb V}_0$ is also closed and
bounded, i.e., compact. Property~\eqref{b13} is checked.

Let us prove property~\eqref{b14}; moreover, property~\eqref{b13a} will be proved
simultaneously. Let us take points $x=(x_1,x_2,\ldots,x_m)$ and $y=(y_1,y_2,\ldots,y_m)
\in {\mathbb V}_0$, i.e., points with coordinates $ x_j,y_j \in (0,2\pi), \ 1\le j\le m$,
satisfying the conditions
\begin{equation}\label{b17}
\prod\limits_{j=1}^m s(x_j) \geq h, \qquad \prod\limits_{j=1}^m s(y_j) \geq h.
\end{equation}
Let us prove that $\dfrac{x+y}{2} \in {\mathbb V}_0$; i.e., $\prod_{j=1}^m s \left(
\dfrac{x_j+y_j}{2}\right) \geq h$. We will prove the more strong inequality
\begin{equation}\label{b18}
\prod\limits_{j=1}^m s\left( \frac{x_j+y_j}{2} \right) \geq \sqrt { \prod\limits_{j=1}^m
s( x_j) \prod\limits_{j=1}^m s( y_j) }\ .
\end{equation}
To do this, let us find the logarithms of the left- and right-hand sides of~\eqref{b18};
we obtain the equivalent inequality
\begin{equation}\label{b20}
\sum\limits_{j=1}^m \ln s\left( \frac{x_j+y_j}{2} \right) \geq \sum\limits_{j=1}^m \frac{
\ln s(x_j)+\ln s(y_j) }{2}.
\end{equation}
The function $\ln s(t)=\ln \left(2 \sin \dfrac{t}{2}\right)$ is strictly convex upwards
on $(0,2\pi)$; therefore, the following inequalities are valid:
\begin{equation}\label{b21}
\ln s\left( \frac{x_j+y_j}{2} \right) \geq \frac{ \ln s(x_j)+\ln s(y_j) }{2},\qquad
j=1,2,\dots,m.
\end{equation}
Consequently, inequality~\eqref{b20} and so inequality~\eqref{b18} are valid. Thus, we
have proved that if $x,y \in {\mathbb V}_0$, then $\dfrac{x+y}{2}\in {\mathbb V}_0$.
Since ${\mathbb V}_0$ is closed, we can conclude that the set ${\mathbb V}_0$ is convex.

Let us ascertain that, in fact, the set ${\mathbb V}_0$ is strictly convex. Let us prove
that if $x,y \in {\mathbb V}_0$ and $x\ne y$, then the following strict inequality holds:
\begin{equation}\label{b22}
{w}\left( \frac{x+y}{2} \right)=\prod\limits_{j=1}^m s \left( \frac{x_j+y_j}{2} \right)
>h.
\end{equation}
If $x\ne y$, then $x_j\ne y_j$ at least for one index~$j$. By the strict convexity of the
function $\ln s(t)=\ln \left(2\sin \dfrac{t}{2}\right)$ on the interval $(0,2\pi)$,
corresponding inequality~\eqref{b21} is strict; but then, inequality~\eqref{b18} is also
strict; consequently, \eqref{b22} holds. The function~$w$ defined by~\eqref{fun-w} is
continuous everywhere in $\mathbb{R}^m$; therefore, there exists a neighborhood
$\mathcal{O}$ of the point $(x+y)/2$ (situated in $(0,2\pi)^m$) in which the inequality
$w(z)>h,\ z\in \mathcal{O},$ is valid. Consequently, this neighborhood lies in the
set~${\mathbb V}_0(h)$.

Since, for $0\le h<2^m$, the set ${\mathbb V}_0(h)$ consists of more than one point, the
previous reasonings imply that ${\mathbb V}_0(h),\ 0<h<2^m$, is strictly convex; in
particular, its interior is nonempty; i.e., this set is a body. Properties~\eqref{b14}
and~\eqref{b13a} are checked. The proof of Lemma~\ref{lemV} is complete.

\ \

{\bf \thesection.3. The intersection of a line with cubes.~ } We assign to a point
$k=(k_1,\ldots,k_m)\in {{\mathbb Z}^m}$ the $m$-dimensional cube
$\mathbb{U}_k=\mathop{\times}\limits_{i=1}^{m} [2\pi k_i,2\pi (k_i+1) )$ in the space
${{\mathbb R}^m},\ m\ge 2$. Evidently, sets from the family $\{{\mathbb{U}_k:\ k \in
{{\mathbb Z}^m}}\}$ are disjoint and $ {{\mathbb R}^m}=\cup\, \{\mathbb{U}_k:\ k \in
{{\mathbb Z}^m}\}. $ Let ${{\mathbb Z}_{0}^m}=\{0\} \times {{\mathbb Z}^{m-1}}$ be the
set of points $k=(0,k_2,\ldots,k_m)\in {{\mathbb Z}^m}$ with integer coordinates the
first of which is zero. For $x\in {{\mathbb R}^m}$, we denote by ${\ell}(x)$ the line
with the directing vector $\mathcal{E}=(1,\ldots,1)$ (i.e., the line orthogonal to the
plane $\mathbb{H}=\mathbb{H}_m$) passing through the point~$x$. In the following lemma,
we study the intersection of lines ${\ell}(x)$ with cubes $\mathbb{U}_k,\ k \in {{\mathbb
Z}_{0}^m}$. We denote by $e_1,e_2,\ldots,e_m$ the unit orts of the space~${{\mathbb
R}^m}. $

\begin{lem}\label{lemCub}
The positional relationship of lines ${\ell}={\ell}(x),\ x\in {{\mathbb R}^m}$, and cubes
${\mathbb{U}_k},\ k \in {{\mathbb Z}_{0}^m}$, has the following properties.

$(1)$ Any line ${\ell}$ intersect at most~$m$ cubes
$$
\{\mathbb{U}_{k^{(p)}}, \ k^{(p)} \in {{\mathbb Z}_{0}^m}\}_{p=1}^{\overline p},\qquad
1\le \overline p\le m.
$$

$(2)$ For $\overline p\ge 2$ and points $k^{(p)} \in {{\mathbb Z}_{0}^m},\ 2\le p\le
\overline p$, the following recurrent formula is valid$:$
\begin{equation} \label{k-p}
k^{(p)}=k^{(p-1)}+\sum\left\{e_i:\ i\in I(p)\right\},\qquad 2\le p\le \overline p,
\end{equation}
where $k^{(1)} \in {{\mathbb Z}_{0}^m}$ and $\{{I}(p)\}_{p=2}^{\overline p}$ is a family
of nonintersecting subsets that form a decomposition of the set of natural numbers
$\{2,\ldots,m\}$.

$(3)$ If a line ${\ell}$ intersect exactly~$m$ cubes $\{\mathbb{U}_{k^{(p)}}, \ k^{(p)}
\in {{\mathbb Z}_{0}^m}\}_{p=1}^{m}$, then
$$
k^{(p)}=k^{(1)}+\sum\limits_{q=2}^{p}e_{i(q)},\qquad 2\le p\le m,
$$
where $\{i(q)\}_{q=2}^m$ is a permutation of the set $\{2,\ldots,m\}$.
\end{lem}

\proof We can write the line ${\ell}={\ell}(x)$ parallel with the vector
$\mathcal{E}=(1,\ldots,1)$ and passing through the point $x=(x_1,\ldots,x_m)$ in the
parametric form
\begin{equation}\label{l(x)}
\lambda(t)=(t,t+x_2-x_1,\ldots,t+x_m-x_1)=t\mathcal{E}+{A},\qquad t\in{\mathbb R},
\end{equation}
where ${A}=\lambda(0)=(0,x_2-x_1,\ldots,x_m-x_1)$; in particular, $\lambda(x_1)=x$. We
denote by $\Pi=\Pi(m)$ the band of points ${z}=(z_1,\ldots,z_m)\in{{\mathbb R}^m}$
satisfying the condition $0\le z_1<2\pi$; we have $\Pi=\cup_{k\in {{\mathbb Z}_{0}^m}}
\mathbb{U}_{k}$. The line ${\ell}$ intersect the boundary hyperplanes $z_1=0$ and
$z_1=2\pi$ of the band $\Pi$ at the points ${A}$ and
$B=\lambda(2\pi)=(2\pi,2\pi+x_2-x_1,\ldots,2\pi+x_m-x_1)$; the distance between these
points is equal to $\sqrt{m}\ 2\pi$. Moreover, the line ${\ell}$ intersect the band $\Pi$
exactly by the half-interval $[{A},B)$. Now we have to find cubes that are intersect this
half-interval.

The fact that the point $\lambda(t)=(t,t+x_2-x_1,\ldots,t+x_m-x_1)$ belongs to the cube
$\mathbb{U}_k,\ k=(0,k_2,\ldots,k_m)\in {{\mathbb Z}_{0}^m}$, means that the
following~$m$ relations are valid:
\begin{equation}\label{l(x)-cub}
0\le t<2\pi; \qquad 2\pi k_i\le t+x_i-x_1<2\pi (k_i+1),\qquad 2\le i\le m.
\end{equation}
First, we construct a cube $\mathbb{U}_{k^{(1)}},\ k^{(1)}\in {{\mathbb Z}_{0}^m}$,
containing the point ${A}=\lambda(0)$; since the cubes $\{{\mathbb{U}_k}\}$ do not
intersect, such cube is unique. Let us choose integers $k^{(1)}_i, \ 2\le i\le m$,
satisfying the condition $2\pi k^{(1)}_i\le x_i-x_1<2\pi (k^{(1)}_i+1),$\ $2\le i\le m$,
and let us set $k^{(1)}=(k^{(1)}_1,k^{(1)}_2,\ldots,k^{(1)}_m)$, where $k^{(1)}_1=0$. It
is clear that $A\in \mathbb{U}_{k^{(1)}}$. Let us consider the numbers $r_i=2\pi
(k_i^{(1)}+1)-(x_i-x_1),\ 1\le i\le m$; they satisfy the condition $0<r_i\le 2\pi$;
moreover, $r_1=2\pi$. Let $T'$ be the set of distinct numbers $\{r_i\}_{i=1}^m$; we
denote by $\overline p,\ 1\le \overline p\le m$, the number of elements of this set. Let
us arrange elements of the set $T'$ in order of magnitude and let us number them by index
$p,\ 2\le p\le \overline p+1$; as a result, we obtain the representation
$T'=\{t_p\}_{p=2}^{\overline p+1}$. In this case, we have $t_{\overline p+1}=2\pi$. Let
us specify the set $T=\{t_p\}_{p=1}^{\overline p}$, where $t_1=0$; elements of this set
are distinct and situated on the half-interval $[0,2\pi)$. We set $A_p=\lambda(t_p),\
1\le p\le \overline p+1$. We have $A_1=\lambda(0)=A$ and $A_{\overline
p+1}=\lambda(2\pi)=B$. For $1\le p\le \overline p$, points~$A_p$ lie in the band~$\Pi$.

First, we assume that $\overline p>1$. For every $p,\ 2\le p\le {\overline p}+1$, we
denote by $I(p)$ the set of (all) indices~$i$ with the property~$r_i=t_p$. The sets
$I(p)$ are disjoint; their union composes the set $\{1,\ldots,m\}$ of the first~$m$
natural numbers. With the help of recurrent relation~\eqref{k-p}, we define a family of
$\overline p-1$ integral points; these points can be also written in the form
$$
k^{(p)}=k^{(1)}+\sum\left\{e_i:\ i\in {\mathcal{I}}(p)\right\},\qquad
{\mathcal{I}}(p)=\bigcup_{q=2}^p I(q),\qquad 2\le p\le \overline p.
$$
For $2\le p\le \overline p$, the sets $I(p)$ do not contain the number 1; therefore, the
points $ k^{(p)},\ 1\le p\le \overline p$, constructed belong to the set~${{\mathbb
Z}_{0}^m}$.

Let us ascertain that the line ${\ell}$ intersect exactly the cubes
$\mathbb{U}_{k^{(p)}},\ 1\le p\le \overline p$. Let us verify that, in fact, the
following stronger assertion is valid:
\begin{equation} \label{l-k-p}
[A_p,A_{p+1})={\ell}\cap \mathbb{U}_{k^{(p)}},\qquad 1\le p\le \overline p.
\end{equation}
On the basis of the choice of values of the parameter $\{t_p\}_{p=1}^{\overline p+1}$,
the definitions of the points $\{A_p\}_{p=1}^{\overline p}$, and the integral points
$\{k^{(p)}\}_{p=1}^{\overline p}$, it is easily seen that the following embeddings are
valid:
\begin{equation} \label{l-k-p-1}
[A_p,A_{p+1})\subset \mathbb{U}_{k^{(p)}},\qquad 1\le p\le \overline p.
\end{equation}
Since, in addition,
$$
\bigcup_{p=1}^{\overline p}[A_p,A_{p+1})=[A,B)={\ell}\cap\Pi,
$$
all embeddings~\eqref{l-k-p-1} turn into equalities~\eqref{l-k-p}. Thus, in fact, the
line~${\ell}$ intersects only the cubes $\mathbb{U}_{k^{(p)}},\ 1\le p\le \overline p$.

If $\overline p=1$, then, as is easily seen, the line ${\ell}$ intersect only one
cube~$\mathbb{U}_{k^{(1)}}$.

Let us discuss the situation when the line ${\ell}$ intersect exactly~$m$ cubes; i.e.,
${\overline p}=m$. This will be in the case if, for any $p,\ 2\le p\le m+1$, the sets
$I(p)$ consist only of one number~$i(p)$; as a consequence, the numbers
$\{i(p)\}_{p=1}^{m}$ form the set $\{1,\ldots,m\}$ and $i(m+1)=1$. The lemma~is proved.

\ \

{\bf \thesection.4. Decomposition of the measure into sum of strictly convex functions.}
By assertion~3 of Lemma~\ref{lem2}, the measure of the set
$[x+a\mathcal{E},x+(a+2\pi)\mathcal{E}) \cap{\mathbb{V}}(h),\ 0<h<2, $ is independent of
the numbers $a\in {\mathbb R}$. In the sequel, it is convenient for us to take $a=-x_1$
for every point $x=(x_1,\dots,x_m)\in \mathbb{H}$. We set ${\mathrm
I}(x)=[x-x_1\mathcal{E},x+(-x_1+2\pi)\mathcal{E})$. Let us consider the function
\begin{equation}
\label{fun-f} {\tau}(x)={\tau}(x;h)={\mathrm{mes}}_1\{[x-x_1\mathcal{E},x+
(-x_1+2\pi)\mathcal{E}) \cap{\mathbb{V}}(h) \}={\mathrm{mes}}_1\{{\mathrm I}(x)
\cap{\mathbb{V}} \}
\end{equation}
of variable $x\in \mathbb{H}$. With this notation, we can write problem~\eqref{b3} in the
form
\begin{equation}\label{fmain}
\Delta_m(h)=\min\{ {\tau}(x;h):\ x \in \mathbb{H} \}.
\end{equation}
Let us decompose function~\eqref{fun-f} into sum of simpler functions. To this end, let
us represent the set ${\mathrm I}(x) \cap{\mathbb{V}}$ in a special form. The
half-interval ${\mathrm I}(x)=[x-x_1\mathcal{E},x+(-x_1+2\pi)\mathcal{E})$ is the
intersection ${\mathrm I}(x)=\ell(x)\cap\Pi$ of the line~$\ell(x)$ and the band
$\Pi=\{z=(z_1,\dots,z_m)\in \mathbb{R}^m:\ 0\le z_1<2\pi\}=[0,2\pi)\times {\mathbb
R}^{m-1}$. Now, using property~\eqref{b12} of the set ${\mathbb{V}}$, we represent the
set ${\mathrm I}(x) \cap{\mathbb{V}}$ as the union of disjoint subsets:
$${\mathrm I}(x) \cap{\mathbb{V}}=
\ell(x)\cap\Pi \cap \bigcup\limits_{k \in {{\mathbb Z}^m}} {\mathbb{V}}_k=
\ell(x) \cap \bigcup\limits_{k \in {{\mathbb Z}^m}}({\mathbb{V}}_k \cap \Pi)=
\ell(x) \cap \bigcup\limits_{k \in {{\mathbb Z}_{0}^m}}{\mathbb{V}}_k=
\bigcup\limits_{k \in {{\mathbb Z}_{0}^m}}(\ell(x) \cap{\mathbb{V}}_k).
$$
We recall that ${{\mathbb Z}_{0}^m}=\{0\} \times {\mathbb Z}^{m-1}$. Thus, we obtain the
representation
\begin{equation}
\label{diff} {\mathrm I}(x) \cap{\mathbb{V}}=\bigcup\left\{\ell(x) \cap {\mathbb{V}}_k:\
{k \in {{\mathbb Z}_{0}^m}}\right\}.
\end{equation}

The sets ${\mathbb{V}}_k$ are strictly convex (see property~\eqref{b14}) and pairwise
disjoint (see property~\eqref{b11});  $\ell(x)$ is a line. Consequently, each of the sets
$\{\ell(x) \cap {\mathbb{V}}_k\}$ can be either the empty set or a point, or a segment;
in addition, the sets $\{\ell(x) \cap {\mathbb{V}}_k\}$ are pairwise disjoint. Therefore,
the following decomposition is valid:
\begin{equation}\label{razl}
{\tau}(x)=\sum\limits_{k \in {{\mathbb Z}_{0}^m}} {\tau}_k(x),\qquad \text{where}\qquad
{\tau}_k(x)={\tau}_k(x;h)={\mathrm{mes}}_1\{\ell(x) \cap {\mathbb{V}}_k \}.
\end{equation}
Let ${\mathbb K}(x)$ be the set of those $k \in {{\mathbb Z}_{0}^m}$ for which $
\ell(x) \cap {\mathbb{V}}_k \ne \varnothing$. By Lemma~\ref{lemCub}, for any point
$x\in {\mathbb{H}}$, the set ${\mathbb K}(x)$ consists of at most~$m$ elements.
Consequently, for any $x\in {\mathbb{H}}$, at most~$m$ summands are nonzero in
sum~\eqref{razl}. In the following lemma, we study properties of the functions
${\tau}_k$, in particular, properties of their supports ${\mathrm{supp}\,}
{\tau}_k\subset \mathbb{H}$.

\begin{lem}\label{prop-fk}
For $m\ge 2,$~ $0<h<2^m$, and $k\in {{\mathbb Z}_{0}^m}$, the functions ${\tau}_k(x),\ x
\in \mathbb{H}$, have the following properties$:$
\begin{equation}\label{vipf}
{\mathrm{supp}\,} {\tau}_k~ \text{is a strictly convex compact body in the hyperplane}~
\mathbb{H};
\end{equation}
\begin{equation}\label{vipsupf1}
{\tau}_k~ \text{is a strictly convex {\rm(}upwards{\rm)} nonnegative function on }~
{\mathrm{supp}\,} {\tau}_k;
\end{equation}
\begin{equation}\label{vipsupf2}
{\tau}_k~ \text{is a continuous function on the hyperplane}~ \mathbb{H}.
\end{equation}
\end{lem}

\proof For a subset $X\in \mathbb{R}^m$, we denote by $\mbox{Pr}_\mathbb{H}(X)$ the
orthogonal projection of the set~$X$ to the hyperplane $\mathbb{H}=\mathbb{H}_m$; i.e.,
$$
\mbox{Pr}_\mathbb{H}(X)=\{x \in \mathbb{H} :\ \ell(x) \cap X \neq \varnothing \}.
$$
In these terms, we have ${\mathrm{supp}\,}
{\tau}_k=\mbox{Pr}_\mathbb{H}({\mathbb{V}}_k)$. By Lemma~\ref{lemV}, the sets
${\mathbb{V}}_k$ are strictly convex and compact bodies in~$\mathbb{R}^m$. Therefore,
each of the sets ${\mathrm{supp}\,} {\tau}_k$ is also a strictly convex and compact
$(m-1)$-dimensional body in the hyperplane~$\mathbb{H}$. The function ${\tau}_k$ is
nonnegative and strictly convex (upwards) on its support ${\mathrm{supp}\,} {\tau}_k$,
since, for $x\in {\mathrm{supp}\,} {\tau}_k$, its value is the measure of the
intersection of the strictly convex set ${\mathbb{V}}_{k}$ and the line~$\ell(x)$. The
function ${\tau}_k$, as a strictly convex function, is continuous in the interior of the
set ${\mathrm{supp}\,} {\tau}_k$. It is easily understood that this function is also
continuous and has zero values at points of the boundary of the set ${\mathrm{supp}\,}
{\tau}_k$. Hence, the function ${\tau}_k$ is continuous on the hyperplane~$\mathbb{H}$.
The lemma~is proved.

\ \

{\bf \thesection.5. The structure of an extremal set.~ }
\begin{lem}\label{tend}
Let $m\ge 2$ and let $0<h<2$. If a point $x^* \in \mathbb{H}$ is extremal
in~$\Delta_m(h)$; more precisely, if an infimum in~\eqref{fmain} is reached at this
point, then the set ${\mathbb K}(x^*)$ consists of~$m$ points and only one of the
functions ${\tau}_k$, $k \in {\mathbb K}(x^*)$, is different from the identical zero.
\end{lem}

\proof Let $L$ be a linear or, more generally, affine subspace of ${{\mathbb R}^m}$ of
dimension $l,\ 1\le l \le m-1$. In the sequel, by the {\it $l$-neighborhood} of a point
$z\in L$ we mean the open ball ${\mathrm O}(z)={\mathrm O}_{r}(z)=\{x\in L:\ |x-z|<{r}\}$
in~$L$ centered at the point~$z$ of certain radius $r>0$; here, $|\cdot|$ is the
Euclidean norm of the space~${{\mathbb R}^m}$.

Let us pay attention to formula~\eqref{razl} for the point $x^*\in \mathbb{H}$ at which
an infimum in~\eqref{fmain} is reached:
\begin{equation}\label{razl-star}
{\tau}(x^*)=\sum\limits_{k \in {\mathbb K}(x^*)} {\tau}_k(x^*).
\end{equation}
We recall that, here, ${\mathbb K}(x^*)$ is the set of those $k \in {{\mathbb Z}_{0}^m}$
for which $ \ell(x^*) \cap {\mathbb{V}}_k \ne \varnothing$; this set contains at most~$m$
elements. The set
$$
{\mathbb{W}}(x^*)=\bigcup \{{\mathbb V}_k:\ k\in {{\mathbb Z}_{0}^m}\setminus{{\mathbb
K}(x^*)}\},
$$
being the union of all the sets ${\mathbb V}_k$ over $ k\in {{\mathbb
Z}_{0}^m}\setminus{{\mathbb K}(x^*)}$, is closed (in ${\mathbb R}^m$) and does not
intersect the line $ \ell(x^*)$. It is easily understood that a distance
$\rho\left(\ell(x^*),\mathbb{W}(x^*)\right)$ between these sets (in the space ${\mathbb
R}^m$) is positive. Let us fix an $(m-1)$-neighborhood ${\mathrm O}_\rho(x^*)=\{x \in
\mathbb{H} :\ |x-x^*|<{\rho}\}$ of the point $x^*$ of radius
$\rho=\rho\left(\ell(x^*),\mathbb{W}(x^*)\right)$ in the hyperplane~$\mathbb{H}$. For
each $ k\in {{\mathbb Z}_{0}^m}\setminus{{\mathbb K}(x^*)}$, the projection
$$
\omega_k=\mbox{Pr}_\mathbb{H}({\mathbb V}_k)={\mathrm{supp}\,} {\tau}_k
$$
of the set ${\mathbb V}_k$ to the hyperplane $\mathbb{H}$ and the neighborhood ${\mathrm
O}_\rho(x^*)$ of the point~$x^*$ are disjoint. Therefore, in fact,
representation~\eqref{razl-star} is valid in the whole neighborhood ${\mathrm
O}_\rho(x^*)$; i.e., the following formula is valid:
\begin{equation}\label{razl-point-x}
{\tau}(x)=\sum\limits_{k\in\, {\mathbb K}(x^*)} {\tau}_k(x),\qquad x\in {\mathrm
O}_\rho(x^*);
\end{equation}
though some of summands in~\eqref{razl-point-x} can be zero.

Points ${k \in {\mathbb K}(x^*)}$ necessarily have one of the following two properties:
(1)~$x^*$ is an interior point of the set $\omega_k={\mathrm{supp}\,} {\tau}_k$;
(2)~$x^*$ is a boundary point of the set $\omega_k={\mathrm{supp}\,} {\tau}_k$. Let us
number the points ${k \in {\mathbb K}(x^*)}$ (or, that is the same, the functions
$\tau_k, {k \in {\mathbb K}(x^*)}$) in a special way. First, note that there exists a
point $k\in {\mathbb K}(x^*)$ with property~(1). Indeed, otherwise,
representation~\eqref{razl-star} implies that ${\tau}(x^*)=0$; consequently,
$\Delta_m({h})=0$. However, this contradicts Lemmas~\ref{lem_delta} and~\ref{lemDelta}.
Thus, let us number points ${k \in {\mathbb K}(x^*)}$ with property~(1) (by upper
indices) from 1 to~$p$. We number points ${k \in {\mathbb K}(x^*)}$ with property~(2), if
they exist, from~$p+1$ to~$q$.

Assume that $m\ge 3$. For every $i=1,\ldots,p$ there exists an $(m-1)$-neighborhood
${\mathrm O}^{(i)}(x^*)$ of the point $x^*$ in the hyperplane $\mathbb{H}$ such that
${\mathrm O}^{(i)}(x^*) \subset {\mathrm{supp}\,} {\tau}_{k^{(i)}}$. For each of the
numbers $i=(p+1),\ldots,q$, in hyperplane~$\mathbb{H}$, there exists an (affine) subspace
$L^{(i)}(x^*)$ (of dimension $m-2$) tangent to ${\mathrm{supp}\,} {\tau}_{k^{(i)}}$ at
the point~$x^*$. Since the set ${\mathrm{supp}\,} {\tau}_{k^{(i)}}$ is strictly convex,
we can assert that this subspace has only one common point with the set
${\mathrm{supp}\,} {\tau}_{k^{(i)}}$, namely, the point~$x^*$. In each of the subspaces
$L^{(i)}(x^*)$, we choose an $(m-2)$-neighborhood ${\mathrm O}^{(i)}(x^*)$ of the
point~$x^*$ (though we can take the subspaces $L^{(i)}(x^*)$ themselves). Note that the
dimension of each of these neighborhoods is equal to~$m-2$. We need to prove that $p=1$
and $q=m$. Let us argue by contradiction. Assume that $p>1$ or $q<m$. We consider the
intersection
$$
\widetilde{{\mathrm O}}(x^*)= \left(\bigcap\limits_{i=1}^q {\mathrm
O}^{(i)}(x^*)\right)\bigcap {\mathrm O}_\rho(x^*)
$$
of the neighborhoods constructed. Let us find the dimension $l=\dim\widetilde{{\mathrm
O}}(x^*)$ of this set. If $p=q$ (i.e., there are no points of the second type),
then~$l=m-1$. Now, let $p<q$. Then, $q-p$ tangent planes intersect at the point~$x^*$ and
$$
l=\dim\widetilde{{\mathrm O}}(x^*)=\dim \bigcap\limits_{i=p+1}^q {\mathrm O}^{(i)}(x^*)
\geq m-1-(q-p)=m-q+p-1.
$$
By assumption, at least one of the inequalities $p>1$ or $q<m$ holds. Therefore, $m-q+p-1
\geq 1$. Thus, $l=\dim \bigcap\limits_{i=p+1}^q {\mathrm O}^{(i)}(x^*)\geq 1$; i.e.,
$\widetilde{{\mathrm O}}(x^*)$ is an $l$-neighborhood of the point~$x^*$ of
dimension~$l\ge 1$. By~\eqref{razl-point-x}, the following representation holds at points
$x \in \widetilde{{\mathrm O}}(x^*)$:
$$
{\tau}(x)=\sum\limits_{i=1}^p {\tau}_{k^{(i)}}(x).
$$
The function ${\tau}$, as a sum of strictly convex functions, is strictly convex on the
set $\widetilde{{\mathrm O}}(x^*)$ and has a minimum on this set (at the point~$x^*$);
the set $\widetilde{{\mathrm O}}(x^*)$ is open in an affine subspace of dimension~$l\ge
1$ (passing through the point~$x^*$). However, this is impossible; a contradiction. In
the case $m\ge 3$, the lemma is proved.

In the case $m=2$, we have to prove that $p=1$ and $q=2$. The negative fact means that
either $p=q=1$ or $p=q=2$. Using the arguments above, we can again ascertain that both
the situations are impossible. Thus, Lemma~\ref{tend} is proved.

As a consequence of Lemma~\ref{tend}, the following assertion is valid which will be
important in the sequel.

\begin{theo}\label{structure}
For $m\ge 1$, $0<h<2$, and the extremal polynomial $P_m\in\mathfrak{P}_m(\Gamma)$ of
problem~\eqref{task2-m-P}, the set $\{ t \in \mathbb{T} :\ |P_m(e^{it})(t)| \ge 1\}$
consists of one segment and~$m-1$ points.
\end{theo}

For $m\ge 2$, this assertion follows from Lemma~\ref{tend}. The case $m=1$ is trivial and
discussed in the first section (see, in particular,~\eqref{delta_h_m1}).

By relation~\eqref{problem-2n-var2} between problems~\eqref{basetproblem} and
~\eqref{task2-m-P}, the following assertion is valid as a special case of
Theorem~\ref{structure}.

\begin{sled} \label{structure-2n}
For $n\ge 1$, $y>1$, and the extremal polynomial $f_n(t)=y \cos nt+f_{n-1}(t)$ of
problem~\eqref{basetproblem}, the set $\{ t \in \mathbb{T} :\ |f_n(t)| \ge 1\}$ consists
of one segment and~$2n-1$ points.
\end{sled}

\section{On Chebyshev polynomials on compact sets of the unit circle.
The completion of studying problem~\eqref{task2-m-P}}

\setcounter{equation}{0}

In this section, we discuss several close problems on polynomials with a fixed leading
coefficient, all zeros of which are situated on the unit circle, and that deviate the
least from zero on compact sets of the circle. Problems of this type are related to
problem~\eqref{task2-m-P}. Using this relationship, Theorem~\ref{structure}, and a result
by L.\,S.\,Maergoiz and
N.\,N.\,Rybakova~\cite{M-R-2008-Conf-Ivanov,M-R-2008-Preprint,M-R-2009-DAN}, we, in
particular, will complete the study of problem~\eqref{task2-m-P}.

\ \

{\bf \thesection.1. Chebyshev polynomials that deviate the least from zero on compact
sets of the unit circle.~ } We recall that, in this paper, we denote by
${\mathfrak{P}}_m(\Gamma)$ the set of algebraic polynomials
\begin{equation} \label{v02-1}
P_m(z)=\prod_{j=1}^m(z-e^{i\phi_j}),\qquad \{\phi_j\}_{j=1}^m\subset \mathbb{R},
\end{equation}
of degree~$m$ with the unit leading coefficient, all~$m$ zeros of which are situated on
the unit circle $\Gamma=\{e^{it}:\, t\in [0,2\pi]\}$ of the complex plane~$\mathbb{C}$.
To the parameter $\theta\in \mathbb{R}$ and polynomial~\eqref{v02-1}, we assign the
polynomial
\begin{equation}
\label{Pm-step} \prod_{j=1}^m(z-e^{i(\phi_j+\theta)})=e^{im\theta}P_m(ze^{-i\theta}),
\end{equation}
whose zeros are obtained by a rotation of zeros of polynomial~\eqref{v02-1} by the
angle~$\theta$ around the origin of the complex plane~$\mathbb{C}$.
Polynomial~\eqref{Pm-step} also belongs to the set $\mathfrak{P}_m(\Gamma)$; we will say
that it is obtained by a rotation of polynomial~\eqref{v02-1} (by the angle~$\theta$
around the origin of the complex plane).

For a compact subset~$Q$ of the circle $\Gamma$, we define the value
\begin{equation} \label{v02-3}
E_m(Q)=\inf\{\|P_m\|_{C(Q)}:\ P_m\in {\mathfrak{P}}_m(\Gamma)\}
\end{equation}
of the best uniform deviation from zero of polynomials from the set
${\mathfrak{P}}_m(\Gamma)$ on~$Q$. For $0<\alpha<\pi$, we denote by
$\mathcal{Q}=\mathcal{Q}(2\alpha)$ the set of all compact sets $Q\subset \Gamma$ whose
(linear) measure $|Q|$ is equal to the number $2\alpha$:\ \ $|Q|=2\alpha$. We are
interested is the least value
\begin{equation} \label{v02-5}
E_m(2\alpha)=E_m(\mathcal{Q}(2\alpha))=\inf\{E_m(Q):\ |Q|=2\alpha\}
\end{equation}
of~\eqref{v02-3} over all compact sets $Q\in \mathcal{Q}(2\alpha)$. The following
assertion  is valid.

\begin{theo}\label{connection}
Problems~\eqref{task2-m-P} and~\eqref{v02-5} are related as follows:
\begin{equation} \label{v02-7}
E_m(2\alpha)=h,\qquad \delta_m(h)=2\pi-2\alpha,\qquad 0<h<2,\qquad 0<\alpha<\pi.
\end{equation}
\end{theo}

This assertion is seemed to be rather natural; let us present some arguments. Let
$Q\subset \Gamma$ be a compact set of measure $|Q|=2\alpha,\ 0<\alpha<\pi$, and let~$P_m$
be a polynomial from ${\mathfrak{P}}_m(\Gamma)$. For $h=\|P_m\|_{C(Q)}$, the measure of
the set $\{t\in\mathbb{T}:\ |P_m(e^{it})|\ge h\}$ is at most~$2\pi-2\alpha$. {\it A
fortiori}, the following inequality is valid:
$$
\delta_m(h)\le 2\pi-2\alpha,\qquad h=\|P_m\|_{C(Q)}.
$$
In the sequel, we will not return to Theorem~\ref{connection} since
Theorems~\ref{probleme-delta-m} and~\ref{Compact-set} proved below contain stronger
assertions in comparison with~\eqref{v02-7}.

Chebyshev polynomials (that deviate the least from zero with the unit leading
coefficient) on compact sets of the complex plane were studied by many mathematicians;
they have numerous applications (see, for example,~\cite{Smirnov}). Let us describe in
more details the results by P.\,L.\,Chebyshev and G.\,Polya on algebraic polynomials that
deviate the least from zero. Let $\mathfrak{P}_m={\mathfrak{P}_m(\mathbb{C})},\, m\ge 1$,
be the set of algebraic polynomials
$$
P_m(x)=x^m+\sum_{k=0}^{m-1}c_kx^k
$$
with the unit leading coefficient and, generally speaking, with complex other
coefficients. P.\,L.\,Chebyshev~\cite{Chebyshev} found the least deviation from zero
\begin{equation} \label{A-Ch-1}
e_m([-1,1])=\inf\{\|P_m\|_{C[-1,1]}:\ P_m\in \mathfrak{P}_m\}
\end{equation}
on the segment $[-1,1]$ of polynomials from the class~${\mathfrak{P}}_m$. Namely, he
showed that
$$
e_m([-1,1])=\frac{1}{2^{m-1}},\qquad m\ge 1,
$$
and the polynomial
$$
P^*_m(x)=\frac{1}{2^{m-1}}\,T_m(x),\qquad T_m(x)=\cos(m\arccos x),\qquad x\in [-1,1].
$$
is extremal. Using linear change of variable, it is easily checked that, for any segment
$I=[a,b]$ of length $|I|=b-a=2{\rho},\ \rho>0$, the quantity
$$
e_m(I)=\inf\{\|P_m\|_{C(I)}:\, P_m\in {\mathfrak{P}_m}\}
$$
has the value $ e_m(I)=2\left(\dfrac{\rho}{2}\right)^m $ and an extremal polynomial can
be found accordingly. For a closed set $Q\subset \mathbb{R}$, we set
\begin{equation}
\label{v02-9} e_m(Q)=\inf\{\|P_m\|_{C(Q)}:\ P_m\in {\mathfrak{P}_m} \}.
\end{equation}
G.\,Polya studied the least value
\begin{equation}
\label{v02-11} e_m(2{\rho})=\inf\{e_m(Q):\ Q\in \mathcal{Q}(2{\rho})\}
\end{equation}
of~\eqref{v02-9} over the family $\mathcal{Q}=\mathcal{Q}(2{\rho})$ of all compact
subsets $Q\subset \mathbb{R}$ of the real line whose measure is equal to a fixed number
$2{\rho},\ \rho>0$. He proved the following assertion (see, for
example,~\cite[p.~23]{Bernshteyn}).

\begin{theo}\label{v02-Polya}
For any $\rho>0$ and any set $Q\in \mathcal{Q}(2{\rho})$, the following inequality is
valid:
$$
e_m(Q)\ge 2\left(\frac{\rho}{2}\right)^m;
$$
it turns into an equality only in the case if~$Q$ is a segment~$($of length~$2{\rho}$$)$.
As a consequence,
$$
e_m(2{\rho})=2\left(\dfrac{\rho}{2}\right)^m.
$$
\end{theo}

Problem~\eqref{v02-5} can be considered to be an analog of problem~\eqref{v02-11}.
However, to study problem~\eqref{v02-5}, we need other arguments in comparison with the
proof of Theorem~\ref{v02-Polya}.

\ \

{\bf \thesection.2. Chebyshev polynomials on an arc of the unit circle.~ } For a segment
$I=[a,b]$ of the real line, we define the arc $\Gamma(I)=e^{iI}=\{e^{it}:\ t\in I\}$ of
length~$|I|$ of the unit circle $\Gamma=\{e^{it}:\ t\in [0,2\pi]\}$. Let
$\mathfrak{P}_m(\Gamma(I))$ be the set of algebraic polynomials~(\ref{v02-1}) with the
unit leading coefficient, all zeros of which are situated on~$\Gamma(I)$. We set
\begin{equation} \label{v02-13}
\varepsilon_m(I)=\min\{\|P_m\|_{C(\Gamma(I))}:\ P_m\in \mathfrak{P}_m(\Gamma(I))\};
\end{equation}
this is one of variants of the problem on polynomials that deviate the least from zero.
Value~\eqref{v02-13} is invariant with respect to shifts of the segment~$I$ (i.e., with
respect to rotations of the arc~$\Gamma(I)$). Therefore, this value depends only on the
length of the segment~$I$. Let us fix the length of segments: $|I|=2{\alpha},\
0<{\alpha}<\pi$, and let us set
\begin{equation} \label{v29-21}
\varepsilon_m(2{\alpha})=\varepsilon_m(I), \qquad |I|=2{\alpha}.
\end{equation}
L.\,S.\,Maergoiz and
N.\,N.\,Rybakova~\cite{M-R-2008-Conf-Ivanov,M-R-2008-Preprint,M-R-2009-DAN} obtained a
solution of problem~\eqref{v29-21}. Earlier, S.\,V.\,Tyshkevich~\cite{Tyshkevich} solved
problem~\eqref{v29-21} for at least two arcs of the circle. His solution is in terms of
the harmonic measure; this solution has a slightly constructive form. A.\,L.\,Lukashov
and S.\,V.\,Tyshkevich also discuss problem~\eqref{v02-13} for several arcs of the circle
in their recent paper~\cite{Lukashov-Tyshkevich-2009}. As a special
case,~\cite{Lukashov-Tyshkevich-2009} contains a solution of problem~\eqref{v29-21} (for
one arc). Note that methods of~\cite{Tyshkevich,Lukashov-Tyshkevich-2009} are different
from that of~\cite{M-R-2008-Conf-Ivanov,M-R-2008-Preprint,M-R-2009-DAN}; in
fact,~\cite{Tyshkevich,Lukashov-Tyshkevich-2009} continue investigations by
A.\,L.\,Lukashov~\cite{Lukashov-2004}. Problem~\eqref{v02-13} on an arc of the circle
without any restrictions on arrangement of zeros was studied earlier in~\cite{Thiran}.
The following assertion is contained
in~\cite{M-R-2008-Conf-Ivanov,M-R-2008-Preprint,M-R-2009-DAN}.

\begin{theo}\label{v30-duga} For $m\ge 1$ and $0<{\alpha}<\pi$, the following formula is
valid:
\begin{equation} \label{v02-15}
\varepsilon_m(2{\alpha})=2\sin^m \frac\alpha2.
\end{equation}
Moreover, with the notation
$$
x_{km}=\cos \frac{\pi(2k-1)}{2m},\qquad a_{k m}(\alpha)\,=1-
2x_{km}^2\cdot\sin^2\frac{\alpha}{2},\qquad k=1, \ldots, n,
$$
for an arc $\Gamma(\alpha)=\{z=e^{it}:\ t\in \mathbb{R},\, |t| \le \alpha \}$ and $m=2n$
or $m=2n+1$, the polynomials
\begin{equation} \label{appr-81}
P_{m}(z)=S_{n}(z),\qquad S_{n}(z)=\prod_{k=1}^{n} (z^2-2a_{k m}(\alpha)\,z+1),\qquad
m=2n,\qquad n\ge 1;
\end{equation}
\begin{equation} \label{appr-82}
P_m(z)=(z-1)S_{n}(z),\qquad m=2n+1,\qquad n\ge 0,
\end{equation}
are the unique extremal polynomials in problem~\eqref{v29-21}.
\end{theo}

For easy references in the sequel, we observe some properties of
polynomials~\eqref{appr-81} and~\eqref{appr-82} on the unit circle. The following two
relations are valid:
\begin{equation} \label{v31-ineq-g-1}
|P_m(e^{it})|\le 2\sin^m \frac{\alpha}{2},\qquad t\in [-\alpha,\alpha],
\end{equation}
\begin{equation} \label{v31-ineq-g-2}
|P_m(e^{it})|>2\sin^m \frac{\alpha}{2},\qquad t\in(\alpha,2\pi-\alpha).
\end{equation}

Starting with polynomials~\eqref{appr-81} and~\eqref{appr-82}, we define the polynomials
\begin{equation} \label{Volga-2r-29a}
g_m(t)=e^{-int}P_m(e^{it}),\qquad m=2n,\qquad n\ge 1;
\end{equation}
\begin{equation} \label{volga-26-2r1-19a}
g_m(t)=-ie^{-i\frac{m}{2}t}P_m(e^{it}),\qquad m=2n+1,\qquad n\ge 0.
\end{equation}
Let us introduce the notation
$$
\lambda=\left(\sin\frac{\alpha}{2}\right)^{-2}.
$$
It is easily checked that polynomial~\eqref{Volga-2r-29a} has the following structure:
\begin{equation} \label{v-2r-extr-g}
g_{2n}(t)=\lambda^{-n} 2T_n(\lambda\cos t-\left(\lambda-1\right)),\qquad t\in\mathbb{R},
\end{equation}
where~$T_n$ is the Chebyshev polynomial (of the first kind) of order~$n$.
Polynomial~\eqref{volga-26-2r1-19a} can be represented in the form
\begin{equation} \label{v-2r1-extr-g}
g_{2n+1}(t))=2\sin \left({\frac{t}{2}}\right)D_n(\lambda\cos t-\left(\lambda
-1\right)),\qquad t\in\mathbb{R},
\end{equation}
where~$D_n$ is an algebraic polynomial of order~$n$ which is expressed in terms of the
Dirichlet kernel on the segment $[-1,1]$ by the formula
$$
D_n(\cos t)=2\mathcal{D}_n(t),\qquad \mathcal{D}_n(t)=\frac{\sin
\left({\dfrac{2n+1}{2}\,t}\right)}{2\sin \left({\dfrac{1}{2}\,t}\right)},\qquad t\ne
2k\pi,\qquad k\in\mathbb{Z}.
$$

Note that polynomial~\eqref{v-2r-extr-g} arose in investigations by
P.\,L.\,Chebyshev~\cite{Chebyshev-III};
before~\cite{M-R-2008-Conf-Ivanov,M-R-2008-Preprint,M-R-2009-DAN}, it was used, in
particular, in papers by A.\,S.\,Mendelev and M.\,S.\,Plotnikov~\cite{b1},~\cite{m}, and
A.\,G.\,Babenko~\cite{BAG}.

\ \

{\bf \thesection.3. Solution of problem~\eqref{task2-m-P}.} The following assertion is
valid for problem~\eqref{task2-m-P}.

\begin{theo}\label{probleme-delta-m}
For $m\ge 1$ and $0<h<2$, the following equality holds:
\begin{equation}\label{z21-b-1}
\delta_m({h})=4\arccos\left(\frac{h}{2}\right)^\frac{1}{m}.
\end{equation}
Moreover, polynomials~\eqref{appr-81} and~\eqref{appr-82} are the unique $($to within an
arbitrary rotation$)$ extremal polynomials in problem~\eqref{task2-m-P} for even and
odd~$m$, respectively; here, the parameters~$h$ and $\alpha$ are related as follows:
\begin{equation}\label{v31-1}
h=2\sin^m \frac{\alpha}{2},\qquad 0<\alpha<\pi.
\end{equation}
\end{theo}

\proof Let
\begin{equation} \label{v30-5}
P_m(z)=\prod_{j=1}^{m}\, \left(z-e^{i\phi_{j}}\right)
\end{equation}
be an extremal polynomial of problem~\eqref{task2-m-P}. By Theorem~\ref{structure}, the
set
$$
\{ t \in \mathbb{T} : |P_m(e^{it})| \ge h\}
$$
consists of $m-1$ points and a segment; we can assume that this segment is symmetrical
with respect to the point $\pi$; more precisely, it has the form $[a,2\pi-a],\ 0<a<\pi$.
On the complementary segment $I^*=[-a,a]$, the inequality $|P_m(e^{it})(t)|\le h, \ t\in
I^*=[-a,a],$ holds or, that is the same, the value of the uniform norm of
polynomial~\eqref{v30-5} on the arc $\Gamma(I^*)=e^{iI^*}=\{e^{it}:\ t\in I^*\}$ is equal
to~$h$:
\begin{equation} \label{v30-7}
\|P_m\|_{C(\Gamma(I^*))}=h.
\end{equation}
Note that, in addition, all~$m$ zeros of polynomial~\eqref{v30-5} belong to the
arc~$\Gamma(I^*)$. By definition~\eqref{v29-21}, assertion~\eqref{v02-15}, and
relation~\eqref{v30-7}, we have
\begin{equation} \label{v30-9}
\varepsilon_m(2{a})=2\sin^m \frac{a}{2}\le h.
\end{equation}
Hence, we obtain the following upper estimate for the length of the segment~$I^*:$
$$
|I^*|=2{a}\le 4\arcsin\left(\frac{h}{2}\right)^\frac{1}{m}.
$$
The equality $\delta_{m}(h)=2\pi-|I^*|$ provides now the estimate
\begin{equation} \label{v30-11}
\delta_{m}(h)=2\pi-2a\ge 2\pi-
4\arcsin\left(\frac{h}{2}\right)^\frac{1}{m}=4\arccos\left(\frac{h}{2}\right)^\frac{1}{m}.
\end{equation}

Using~\eqref{v31-1}, we represent the parameter~$h$ in terms of $\alpha\in(0,\pi)$. By
properties~\eqref{v31-ineq-g-1} and~\eqref{v31-ineq-g-2}, polynomials~\eqref{appr-81}
and~\eqref{appr-82} provide the inverse estimate. Thus, assertion~\eqref{z21-b-1} and the
property of polynomials~\eqref{appr-81} and~\eqref{appr-82} to be extremal are proved. It
is seen from the proof that $a=\alpha$ and polynomial~\eqref{v30-5} solves
problem~\eqref{v29-21}. By Theorem~\ref{v30-duga}, polynomial~\eqref{v30-5} coincides
with~\eqref{appr-81} or~\eqref{appr-82} depending on the evenness of the number~$m$.
Theorem~\ref{probleme-delta-m} is proved.

\ \

{\bf \thesection.4. The investigation of problem~\eqref{v02-5}.~ } Let us return to
approximation problem~\eqref{v02-5}. Evidently, values~\eqref{v02-5} and~\eqref{v29-21}
are related by the inequality
\begin{equation} \label{AQ-p-7}
E_m(2\alpha)\le
\varepsilon_m(\Gamma(\alpha))=\varepsilon_m(2\alpha)=2\sin^m\frac{\alpha}{2},\qquad
0<\alpha<\pi.
\end{equation}
Now, we will see that, in fact, they coincide.

\begin{theo}\label{Compact-set} For any $\alpha,\ 0<\alpha<\pi$, and
any set $Q\in \mathcal{Q}(2{\alpha})$, the following inequality is valid:
$$
E_m(Q)\ge2\sin^m\frac{\alpha}{2};
$$
it turns into an equality only in the case if the set~$Q$ is an arc~$($of length
$2{\alpha})$. As a consequence,
$$
E_m(2{\alpha})=\varepsilon_m(2\alpha)=2\sin^m\frac{\alpha}{2}.
$$
\end{theo}

\proof Assume that, for a compact subset $Q\subset \Gamma$ of measure $|Q|=2\alpha,\
0<\alpha<\pi$, the following inequality is valid:
\begin{equation} \label{25okt09}
E_m(Q)\le h,\quad{\mbox{where}}\quad h=2\sin^m\frac{\alpha}{2};
\end{equation}
by definitions~\eqref{v02-3} and~\eqref{v02-13}, an arbitrary arc of the unit circle of
length~$2\alpha$ has this property  {\it a~fortiori}. Let us prove that, then, the
set~$Q$ is sure an arc of length~$2\alpha$ and inequality~\eqref{25okt09} turns into an
equality. Thus, Theorem~\ref{Compact-set} will be proved.

Let ${P}_m\in \mathfrak{P}_m(\Gamma)$ be the polynomial on which an infimum
in~\eqref{v02-3} is reached for the set~$Q$ under consideration. We use the notation
$h'=\|{P}_m\|_{C(Q)}$; assumption~\eqref{25okt09} means that~$h'\le h$. Let us ascertain
that the following estimate is valid for the measures of the set $ \{z\in \Gamma:\
|{P}_m(z)|\geq h'\}$:
\begin{equation} \label{25okt09-2}
|\{z\in \Gamma:\ |{P}_m(z)|\geq h'\}|\le 2\pi-2\alpha.
\end{equation}
Indeed,
$$
\{z\in \Gamma:\ |{P}_m(z)|\geq h'\}=\{z\in \Gamma:\ |{P}_m(z)|>h'\}\cup \{z\in \Gamma:\
|{P}_m(z)|=h'\}.
$$
The embedding $\{z\in \Gamma:\ |{P}_m(z)|>h'\}\subset \Gamma\setminus Q$ is valid;
therefore,
$$
|\{z\in \Gamma:\ |{P}_m(z)|>h'\}|\le 2\pi-2\alpha.
$$
Using, for instance, representation~\eqref{G-s}--\eqref{bb2} of the polynomial~$P_m$, we
can easily check that the set $\{z\in \Gamma:\ |{P}_m(z)|=h'\}$ is finite and so has the
measure zero. Therefore, estimate~\eqref{25okt09-2} is really valid.

By definition~\eqref{task2-m-P} and inequality~\eqref{25okt09-2}, we have
\begin{equation} \label{25okt09-3}
\delta_m(h')\le |\{z\in \Gamma:\ |{P}_m(z)|\geq h'\}|\le 2\pi-2\alpha.
\end{equation}
Using~\eqref{z21-b-1}, we can easily check that if $h=2\sin^m\alpha/2$, then
$2\pi-2\alpha=\delta_m(h)$. Therefore, inequality~\eqref{25okt09-3} can be written in the
form
\begin{equation} \label{25okt09-4}
\delta_m(h')\le \delta_m(h).
\end{equation}
By Theorem~\ref{probleme-delta-m}, value~\eqref{task2-m-P} decreases with respect to its
argument; therefore,~\eqref{25okt09-4} implies that $h'\ge h$. Taking into account the
property $h'\le h$, we conclude that $h'=h$ and
$$
\delta_m(h)=|\{z\in \Gamma:\ |{P}_m(z)|\geq h\}|.
$$
This fact means that the polynomial ${P}_m$ is extremal in problem~\eqref{task2-m-P}. By
Theorem~\ref{probleme-delta-m}, the polynomial ${P}_m$ coincides to within a rotation
with~\eqref{appr-81} or~\eqref{appr-82} depending on the evenness of the number~$m$.
Consequently, the set
\begin{equation} \label{v03-5}
\{z\in \Gamma:\ |{P}_m(z)|\le h\},\qquad h=2\sin^m\frac{\alpha}{2},
\end{equation}
is an arc of length~$2\alpha$. The set~$Q$ is compact; its measure is also equal
to~$2\alpha$; this set belongs to arc~\eqref{v03-5}. Therefore, the set~$Q$ coincides
with arc~\eqref{v03-5}. The arguments above contain the equality $h=\|P_m\|_{C(Q)}$,
which means that~\eqref{25okt09} turns into an equality for the set or, more precisely,
the arc~$Q$. Theorem~\ref{Compact-set} is proved.

\section{The completion of the proof of Theorem~\ref{maint}.\\
Trigonometric polynomials deviating the least from zero\\
with respect to the uniform norm on compact subsets\\ of the torus that have a given
measure}

\setcounter{equation}{0}

{\bf \thesection.1. Proof of Theorem~\ref{maint}.~ } Assertion~\eqref{Eq} follows from
Corollary~\ref{sigma-delta} and Theorem~\ref{probleme-delta-m}, more precisely, from
equalities~\eqref{problem-2n-var2} and~\eqref{z21-b-1}. It remains to describe the set of
extremal polynomials of problem~\eqref{basetproblem}. By Lemma~\ref{lem1}, an extremal
polynomial has only real roots; i.e., it belongs to the set
$\mathcal{F}_{n}^{\,real}(y)$. According to the results of the first section
(see~\eqref{T-real} and~\eqref{a6-y-refl}), a polynomial $f_{n}\in
\mathcal{F}_{n}^{\,real}(y)$ has the representation
\begin{equation}
\label{P-predstavlenie} f_{n}(t)=\frac{y}{2} e^{-int} P_{2n}(e^{it}),
\end{equation}
where
\begin{equation}
\label{P-predstavlenie-0} P_{2n}(z)=\prod_{j=1}^{2n}(z-e^{i\phi_j})
\end{equation}
is a polynomial from the set $\mathfrak{P}_{2n}(\Gamma)$ with the property
\begin{equation}
\label{P-predstavlenie-1} \Phi=\sum_{j=1}^{2n}\phi_j=2\pi k,\qquad k\in \mathbb{Z}.
\end{equation}
By Corollary~\ref{sigma-delta}, polynomial~\eqref{P-predstavlenie} is extremal in
problem~\eqref{basetproblem} if and only if polynomial~\eqref{P-predstavlenie-0} is
extremal in problem~\eqref{task2-m-P} for $ m=2n$ and $h={2}/{y}$. According to
Theorem~\ref{probleme-delta-m}, such a polynomial coincides to within a rotation with
polynomial~\eqref{appr-81} if  the parameters satisfy relations~\eqref{v31-1}.

Polynomial~\eqref{appr-81} has~$n$ pairs of complex-conjugate roots $e^{\pm i\phi_j}$,
where $0<\phi_j<\alpha<\pi,\ 1\le j\le n$. For this polynomial,
sum~\eqref{P-predstavlenie-1} is equal to zero; i.e., $k=0$. Hence, on the base of
formulas~\eqref{Volga-2r-29a},~\eqref{v-2r-extr-g}, and~\eqref{P-predstavlenie}, we
conclude that the polynomial
\begin{equation}
\label{extr-polinom-BAG-0} f_{n,0}(t)=\frac{y}{2} e^{-int} S_{n}(e^{it})=T_{n}\left(
y^{\frac{1}{n}}\cos t-y^{\frac{1}{n}}+1\right)
\end{equation}
belongs to the set $\mathcal{P}_{2n}(y)$ and is extremal in problem~\eqref{basetproblem}.

The procedure of rotation~\eqref{Pm-step} of polynomial~\eqref{appr-81} by a value
$\theta\in \mathbb{R}$ gives the polynomial
\begin{equation}
\label{extr-polinom-BAG-2} P_{2n}(z)=e^{i2n\theta}S_{n}(ze^{-i\theta}),
\end{equation}
for which sum~\eqref{P-predstavlenie-1} is equal to the number~$2n\theta$.
By~\eqref{P-predstavlenie-1}, for the respective polynomial~\eqref{P-predstavlenie} to be
extremal it is necessary and sufficient to have $2n\theta=2\pi k,\ k\in \mathbb{Z}$, or
$$
\theta=\theta_k=\frac{k\pi}{n},\qquad k\in \mathbb{Z}.
$$
For this value $\theta$, we have
$$
f_n(t)=f_{n,k}(t)=\frac{y}{2} e^{-int} S_{n}(e^{i(t-\theta_k)})=\frac{y}{2}
e^{in\theta_k} e^{-in(t-\theta_k)} S_{n}(e^{i(t-\theta_k)})=(-1)^kf_{n,0}(t-\theta_k)
$$
or, by~\eqref{extr-polinom-BAG-0},
\begin{equation}
\label{extr-polinom-BAG-k} f_n(t)=f_{n,k}(t)=(-1)^kT_{n}\left( y^{\frac{1}{n}}\cos
(t-\theta_k)-y^{\frac{1}{n}}+1\right).
\end{equation}
Thus, we have shown that extremal polynomials of problem~\eqref{basetproblem} are
described by formula~\eqref{extr-polinom-BAG-k}; this is precisely the assertion of
Theorem~\ref{maint}. Theorem~\ref{maint} is proved.

\ \

{\bf \thesection.2. On trigonometric polynomials that deviate the least from zero on
compact sets of a given measure.~ } For $0<\alpha<\pi$, we denote by
$\mathcal{T}(2\alpha)$ the set of all compact subsets~$Q$ of the torus ${\mathbb{T}}$
whose measure $|Q|$ is equal to the number~$2\alpha$:\ \ $|Q|=2\alpha$. For $n\ge 1$ and
a compact subset $Q\subset \mathbb{T}$, we define the value
\begin{equation} \label{v29-11}
U_n(Q)=\inf\{\|\cos nt-f_{n-1}\|_{C(Q)}:\ f_{n-1}\in \mathcal{F}_{n-1}\}
\end{equation}
of the best uniform approximation of the function $\cos nt$ by the
family~$\mathcal{F}_{n-1}$ of trigonometric polynomials of order~$n-1$ on the set~$Q$. We
are interested in the least value
\begin{equation} \label{v29-13}
U_n(2\alpha)=U_n(\mathcal{T}(2\alpha))=\inf\{U_n(Q):\ |Q|=2\alpha\}
\end{equation}
of~\eqref{v29-11} over all compact sets~$Q\in \mathcal{T}(2\alpha)$.

Problems~\eqref{v29-11} and~\eqref{v29-13} can also be considered as analogs of
problems~\eqref{A-Ch-1} and~\eqref{v02-11} studied by P.\,L.\,Chebyshev and G.\,Polya.
However, for the study of~\eqref{v29-11} and~\eqref{v29-13}, other methods are applied.
A.\,L.\,Lukashov~\cite{Lukashov-2004} gave a solution of problem~\eqref{v29-11} for a
finite set of segments; however, terms that he used to obtain these results do not allow
one to conclude anything about problem~\eqref{v29-13}.

Problem~\eqref{v29-11} for the segments
\begin{equation} \label{v08-14}
I_k(2\alpha)
=I(2\alpha)+\frac{k\pi}{n}=\left[-\alpha+\frac{k\pi}{n},\alpha+\frac{k\pi}{n}\right],\qquad
k\in \mathbb{Z},
\end{equation}
that are shifts of the segment $I(2\alpha)=[-\alpha,\alpha]$, plays an important role.
For these segments, a solution of problem~\eqref{v29-11} can be easily given.

\begin{lem}\label{otrezok}
Let $n\ge 1$ and $0<\alpha<\pi$. Then,
\begin{equation}\label{appr-10}
U_n\left(I_k(2\alpha)\right)=\sin^{2n} \frac{\alpha}{2},\qquad k\in \mathbb{Z},
\end{equation}
and the polynomials
\begin{equation}\label{v09-10}
{\widetilde f}_{n,k}(t)=(-1)^k\left(\sin^{2n} \frac{\alpha}{2}\right)\cdot T_{n}\left(
y^{\frac{1}{n}}\cos \left(t+\frac{\pi k}{n}\right)-y^{\frac{1}{n}}+1\right),\qquad
y=\sin^{-2n} \frac{\alpha}{2},
\end{equation}
are extremal; they differ from polynomials~\eqref{extr-polinom-BAG} only by the proper
normalization.
\end{lem}

\proof We restrict our attention to the case~$k=0$. The polynomial
\begin{equation}\label{v09-10a}
{\widetilde f}_{n}(t)={\widetilde f}_{n,0}(t)=\left(\sin^{2n}
\frac{\alpha}{2}\right)\cdot T_{n}\left( y^{\frac{1}{n}}\cos t-y^{\frac{1}{n}}+
1\right),\qquad y=\sin^{-2n} \frac{\alpha}{2},
\end{equation}
has the form ${\widetilde f}_{n}(t)=\cos nt+f_{n-1}(t),\ f_{n-1}\in \mathcal{F}_{n-1}$;
it has a $(2n+1)$-point alternance on the segment $[-\alpha,\alpha]$. Therefore,
$U_n(I_k(2\alpha))=\|{\widetilde f}_{n}\|_{C[-\alpha,\alpha]}=\sin^{2n}{\alpha}/{2}$. The
lemma~is proved.

\ \

The following assertion containing solution of problem~\eqref{v29-13} is valid.

\begin{theo}\label{v29-approx-II}
For any $\alpha,\ 0<\alpha<\pi$, for any compact subset $Q\subset \mathbb{T}$ of the
torus of measure $|Q|=2\alpha$, the following inequality is valid:
\begin{equation} \label{v08-13}
U_n(Q)\ge \sin^{2n} \frac{\alpha}{2};
\end{equation}
it turns into an equality only on segments~\eqref{v08-14}. As a consequence, the
following equality holds for value~\eqref{v29-13}:
\begin{equation} \label{AQ-p-8}
U_n(2\alpha)=\sin^{2n} \frac{\alpha}{2}.
\end{equation}
\end{theo}

\proof The proof of this assertion is carried out with the help of Theorem~\ref{maint} by
the same scheme as the proof of Theorem~\ref{Compact-set} was carried out, starting from
Theorem~\ref{probleme-delta-m}. Indeed, let us assume that, for a compact subset
$Q\subset \mathbb{T}$ of measure $|Q|=2\alpha,\ 0<\alpha<\pi$, the inequality $U_m(Q)\le
\sin^{2n} {\alpha}/{2}$ is valid. Let $f_n(t)=\cos nt-f_{n-1},\ f_{n-1}\in
\mathcal{F}_{n-1}$, be a polynomial on which an infimum in~\eqref{v29-11} is reached for
this set~$Q$. We have $d'=\|f_n\|_{C(Q)}\le d=\sin^{2n}\alpha/2$. Let us estimate the
measure of the set
$$
\{t\in \mathbb{T}:\ |y'\,f_n|\geq 1\},\qquad y'=1/d',
$$
from above. This set can be represented in the form
$$
\{t\in \mathbb{T}:\ |y'\,f_n|\geq 1\}=\{t\in \mathbb{T}:\ |y'\,f_n|>1\}\cup \{t\in
\mathbb{T}:\ |y'\,f_n|=1\}.
$$
The embedding $\{t\in \mathbb{T}:\ |y'\,f_n|>1\}\subset \mathbb{T}\setminus Q$ is valid;
consequently, $|\{t\in \mathbb{T}:\ |y'\,f_n|>1\}|\le 2\pi-2\alpha$. However, the set
$\{t\in \mathbb{T}:\ |y'\,f_n|=1\}$ is finite; so, $|\{t\in \mathbb{T}:\
|y'\,f_n|=1\}|=0$. Therefore, the following estimate is valid:
$$
|\{t\in \mathbb{T}:\ |y'f_n(t)|\geq 1\}|\le 2\pi-2\alpha.
$$

The function $y'f_n$ belongs to the set $\mathcal{F}_n(y')$ and satisfies the following
inequalities:
\begin{equation} \label{v08-15}
\sigma_n(y')\le |\{t\in \mathbb{T}:\ |y'f_n|\geq 1\}|\le 2\pi-2\alpha=\sigma_n(y),\qquad
y=1/d=\sin^{-2n}\frac{\alpha}{2}.
\end{equation}
By Theorem~\ref{maint}, value~\eqref{basetproblem} increases with respect to its
argument. Since $y'\ge y$, \eqref{v08-15} implies that $y'=y$ and
\begin{equation} \label{v08-17}
\sigma_n(y)=|\{t\in \mathbb{T}:\ |yf_n(t)|\geq 1\}|=2\pi-2\alpha.
\end{equation}
This fact means that the polynomial $y f_n$ is extremal in problem~\eqref{basetproblem}.
By Theorem~\ref{maint}, the polynomial $y f_n$ coincides with one of
polynomials~\eqref{extr-polinom-BAG}. Consequently, the set
\begin{equation} \label{v08-19}
\{t\in \mathbb{T}:\ |yf_n(t)|\ge 1\}=\{t\in \mathbb{T}:\ |f_n(t)|\ge d\},\qquad
d=\sin^{2n}\frac{\alpha}{2},
\end{equation}
is one of segments~\eqref{v08-14} of length~$2\alpha$. The set~$Q$ is compact; its
measure is also equal to~$2\alpha$; this set belongs to segment~\eqref{v08-19}.
Consequently,~$Q$ coincides with segment~\eqref{v08-19}; i.e., it coincides with one of
segments~\eqref{v08-14}. In this case, inequality~\eqref{v08-13} turns into an equality.
Theorem~\ref{Compact-set} is proved.

As a consequence of Theorems~\ref{maint} and~\ref{v29-approx-II}, the following analog of
Theorem~\ref{connection} is valid.

\begin{sled}\label{v08-t2}
Problems~\eqref{basetproblem} and~\eqref{v29-13} are related as follows:
\begin{equation} \label{v08-3}
U_n(2\alpha)=y^{-1},\qquad 2\alpha=2\pi-\sigma_n\left(y\right),\qquad y>1.
\end{equation}
\end{sled}

{\bf \thesection.3. Sharp constant in inequality~\eqref{babenko-ineq}.~ }
Theorem~\ref{maint} allow us to find a value of the best constant~$\beta_n$ in
inequality~\eqref{babenko-ineq}.

\begin{theo}\label{Ineq-Babenko}
For any $n\ge 1$, the following formula is valid for the best constant~$\beta_n$ in
inequality~\eqref{babenko-ineq}:
\begin{equation}
\label{AM-p-5}
\beta_{n}=\sqrt{2n}.\end{equation}
\end{theo}

\proof The leading harmonic of trigonometric polynomial~\eqref{trig_polin} can be written
in the form $a_{n}\cos nt+b_{n}\sin nt=y\cos (nt+t_n)$, where
$y=y(f_n)=\sqrt{a_n^2+b_n^2,}$ and~$t_n$ is the respective shift of the argument.
Functional~\eqref{mera} is invariant with respect to a shift of argument of the
polynomial; hence, we can assume that $y\cos nt$ is the leading harmonic of the
polynomial~$f_n$; i.e., $f_n\in {\mathcal{F}_n(y)}$. Studying
inequality~\eqref{babenko-ineq}, we have to restrict our attention only to polynomials
with~$y=y(f_n)>1$; in addition, it is reasonable to choose lower harmonics of the
polynomial so that functional~\eqref{mera} have the least value. As a result, we arrive
at the following representation of the best constant~$\beta_n$ in
inequality~\eqref{babenko-ineq}:
$$
\beta_{n}=\sup\limits_{y>1} \frac{\mu(y\cos(nt))} {\sigma_n(y)}.
$$
We have
$$
\mu(y\cos(nt))=4\arccos \frac{1}{y}.
$$
Now, applying Theorem~\ref{maint}, we obtain
$$
\beta_{n}=
\sup\limits_{y>1} \frac{4\arccos \frac{1}{y}} {4\arccos\frac{1}{y^{\frac{1}{2n}}}}=
\sup\limits_{0\le t<1} \frac{\arccos {t^{2n}}} {\arccos{t}}=\lim_{t\to 1-0} \frac{\arccos
{t^{2n}}} {\arccos{t}}=\sqrt{2n}.
$$
Assertion~\eqref{AM-p-5} is proved.

\ \

\noindent{\bf Remark.} Problems~\eqref{basetask},~\eqref{v29-11}, and~\eqref{v29-13} for
leading harmonic of the general form $A\cos nt+B\sin nt$ are reduced to the case $y\cos
nt$ considered in this paper by an appropriate change of variable.

\section*{Acknowledgments}

The authors are grateful to R.R.~Akopyan, A.G.~Ba\-benko, and
P.Yu.~Gla\-zy\-ri\-na for careful reading of the manuscript and
useful discussions. This work was supported by the Russian
Foundation for Basic Research (project No. 08-01-00213) and by
the Program for State Support of Leading Scientific Schools of
the Russian Federation (project No. NSh-1071.2008.1).

\end{document}